\documentclass[11pt,a4paper]{article}

\usepackage{amsmath}
\usepackage{amssymb}
\usepackage{amsfonts}

\usepackage[utf8]{inputenc}

\usepackage{graphicx}
\usepackage[autostyle]{csquotes}
\usepackage[linktocpage]{hyperref}
\usepackage{fullpage}
\usepackage{slashed}
\usepackage{mathtools}
\usepackage{indentfirst}
\usepackage{bm}
\usepackage{mathrsfs} 
\usepackage{authblk}
\usepackage{tikz}
\usepackage{float}

\allowdisplaybreaks

\newtheorem{theorem}{Theorem}[section]

\newtheorem{Defn}[theorem]{Definition}
\newtheorem{Thm}[theorem]{Theorem}
\newtheorem{Prop}[theorem]{Proposition}
\newtheorem{Lem}[theorem]{Lemma}

\newtheorem{Rmk}[theorem]{Remark}
\newtheorem{Conj}[theorem]{Conjecture}

\newcommand{\intg}{\mathbb{Z}}
\newcommand{\rational}{\mathbb{Q}}
\newcommand{\real}{\mathbb{R}}
\newcommand{\complex}{\mathbb{C}}

\newcommand{\lb}{\left(}
\newcommand{\rb}{\right)}
\newcommand{\lsb}{\left[}
\newcommand{\rsb}{\right]}
\newcommand{\lac}{\left\{}
\newcommand{\rac}{\right\}}

\newcommand{\ptl}{\partial}

\newcommand{\RN}[1]{%
  \textup{\uppercase\expandafter{\romannumeral#1}}}

	\newcommand{\rotsimeq}{\rotatebox[origin=c]{90}{$\simeq$}}

\title{A supergroup series for knot complements}
\author{John Chae}
\affil{yjchae@formerstudents.ucdavis.edu}
\date{}  

\begin{document}

\maketitle

\begin{abstract}

We introduce a three variable series invariant $F_K (y,z,q)$ for plumbed knot complements associated with a Lie superalgebra $sl(2|1)$. The invariant is a generalization of the $sl(2|1)$-series invariant $\hat{Z}(q)$ for closed 3-manifolds introduced by Ferrari and Putrov and an extension of the two variable series invariant defined by Gukov and Manolescu (GM) to the Lie superalgebra. We derive a surgery formula relating $F_K (y,z,q)$ to $\hat{Z}(q)$ invariant. We find appropriate expansion chambers for certain infinite families of torus knots and compute explicit examples. Furthermore, we provide evidence for a non semisimple $Spin^c$ decorated TQFT from the three variable series. We observe that the super $F_K (y,z,q)$ itself and its results exhibit distinctive features compared to the GM series.

\end{abstract}

\tableofcontents

\section{Introduction}

Topological quantum field theories (TQFTs) have been a fruitful source of the interactions between physics and topology. From one to four dimensions, TQFTs have provided physical realizations of topological invariants or predicted new ones. Examples include colored Jones polynomials, HOMFLY-PT polynomials of links~\cite{W1,NRZ}, Donaldson invariants and Seiberg-Witten invariants of smooth four manifolds~\cite{W3,W4}. In three dimensions, Chern-Simons TQFT predicted the Witten-Reshetikhin-Tureav (WRT)-invariant of 3-manifolds~\cite{W1}. The introduction of this invariant motivated a rigorous construction of the invariant via quantum group $U_q (sl(2))$ and their representations~\cite{RT}. This in turn has led to the quantum R-matrix method for computations of the link polynomials. 
\newline

On the mathematics side, TQFT was axiomatized in \cite{At,Se} (see \cite{F} for a review) and its breadth and depth have been enriched. One direction of advancement of TQFTs has been constructions of extended TQFTs. There has been progress in the classification of such TQFTs~\cite{BD,L}. Another line of development of TQFTs is constructions of non semisimple TQFTs associated with a variety of quantum groups. In three dimensions, this kind of TQFTs used non semisimple categories~\cite{GKP2} and the modified quantum dimension~\cite{GPT,GKP1}. A non semisimple TQFT has produced a new non semisimple quantum invariant of links and 3-manifolds called CGP invariant~\cite{CGP}. Advantages of the non semisimple invariants are that they can distinguish manifolds that are not feasible by semisimple invariants and they yield nonzero results in cases the latter vanish. The underlying quantum groups of the TQFTs have been generalized to quantum supergroups~\cite{GP1,GP2,H}.
\newline

Another rich source of the interactions between physics and topology is the categorification program~\cite{CF} (see \cite{A,S} for reviews). It has not only deepened understanding of quantum invariants of manifolds but it also provided powerful tools. In case of link polynomials, they were turned out to be graded Euler characteristics of homology theories. For example, Jones polynomials and HOMFLY-PT polynomials are Euler characteristics of Khovanov (co)homology~\cite{K1, K2} and Khovanov-Rozansky homology~\cite{KR2}, respectively. Furthermore, quantum group itself was categorified, which combined with quantum Weyl group have led to a different approach for computing link polynomials~\cite{KL2,La}.
\newline

From the physics perspective of categorification, string theory has played a vital role (see \cite{G} for a review). In the case of Khovanov homology, a brane system from string/M theory was constructed in \cite{W2}. A physical realization of Khovanov-Rozansky homology was achieved through an application of topological string theory~\cite{GSV}.
\newline

A major challenge of the categorification program has been categorifying the WRT invariant of closed 3-manifolds $Y$. The invariant is defined at root of unity and does not have manifest integrality property to be the Euler characteristic of a homology theory. A strategy for categorification has been proposed in \cite{K3,EQ}. On the physics side, a 3-dimensional supersymmetric QFT originated from 6 dimensions predicted an existence of a power series with integral coefficients associated with the WRT invariant~\cite{GPV,GPPV}. This $q$ power series was denoted by $\hat{Z}_b $ and labeled by $Spin^c$ structures of $Y$. In addition, $\hat{Z}_b$ is associated with a Lie algebra $sl(2)$ and itself is a topological invariant of $Y$, which is a vast generalization of \cite{LZ}. It was conjectured that $\hat{Z}_b $ decomposes the WRT invariant as a linear combination. This was proven for a particular class of 3-manifolds~\cite{Mu}. Importantly, it was conjectured that $\hat{Z}_b$ is the graded Euler characteristic of a homology theory that provides the desired categorification of the WRT invariant. 
\newline

A generalization to 3-manifolds with torus boundary, in particular, plumbed knot complements, was achieved in \cite{GM}. This resulted in a two variable series invariant $F_K(x,q)$ for a complement of a knot $K$. There have been extensive developments in both $F_K$ and $\hat{Z}$. For example, extensions to higher rank Lie groups~\cite{P1}, R-matrix and state sum approach~\cite{P2,P3,Gr}, satellite knots~\cite{C3}, and quantum modularity property~\cite{CCFGM, CCKPS} (see also \cite{GHNPPS,EGGKPS,GPP,C2,C4,GJ} and references therein).
\newline

Motivated by $\hat{Z}_b$, its extension to a Lie superalgebra was introduced in \cite{FP}. In case of $sl(2|1)$, a new $q$ power series was introduced and was denoted by $\hat{Z}_{b,c} (q)$ and carries two labels $(b,c) \in Spin^c (Y) \times Spin^c (Y)$. For a class of 3-manifolds called plumbed manifolds $Y(\Gamma)$ , it was shown that $\hat{Z}_{b,c}$ decomposes a quantum invariant of $Y(\Gamma)$ constructed in \cite{H} (see Section 2 for a review). From physics viewpoint, string/M theory predicted that the existence of $\hat{Z}_{b,c} (q)$ and it was claimed to be a topological invariant of $Y(\Gamma)$ (see Appendix C for details).
\newline

In this paper, we generalize $\hat{Z}_{b,c} (q)$ to a complement of $K$ motivated by \cite{GM}. In particular, we introduce a three variable power series invariant super $F_K (y,z,q)$ for plumbed knot complements, derive a surgery formula that allows to connect to $\hat{Z}_{b,c} (q)$ and compute examples for torus knots. We will observe that super $F_K (y,z,q)$ is qualitatively different from $F_K(x,q)$ associated with $sl(2)$. Furthermore, we show that $\hat{Z}_{b,c} (q)$ is a topological invariant of $Y(\Gamma)$. 
\newline

\noindent \textbf{Statement of Results} We begin with plumbed knot complements $Y_K$ that are represented by plumbing graphs with one distinguished vertex. For plumbing graphs satisfying (weakly) negative definite condition, we obtain an invariant 
$$
\hat{Z}_{b,c}(Y_K ; y,z,n,m,q;\alpha_i).
$$
This is a series in three variables $y,z,q$ and depends on the choice of relative $Spin^c$ structures $(b,c) \in Spin^c (Y_K,\ptl Y_K) \times Spin^c (Y_K,\ptl Y_K)$ and of chambers $\alpha_i, i=\pm$. Furthermore, it also depends on $n,m \in \intg$. Under gluing of the knot complements, the series behaves as follows.
\begin{Thm} Let $Y_1$ and $Y_2$ be knot complements represented by (weakly) negative definite plumbing graphs and $Y=Y_1 \cup_{T^2} Y_2$ be the result of gluing them along their common torus boundary. Let also $(b_1 , c_1)$ and $(b_2 , c_2)$ be relative $Spin^c$ structures of $Y_1$ and $Y_2$, respectively, which results in $Spin^c$ structures $(b,c)$ of $Y$. The gluing yields
$$
\hat{Z}_{b,c}[Y;q] = (-1)^{\tau} q^{\chi} \sum_{n,m} \int \frac{dy}{i2\pi y}\frac{dz}{i2\pi z}\, \hat{Z}^{(\alpha_i)}_{b_1 ,c_1} (Y_1 ; y,z,n,m,q)\, \hat{Z}^{(\alpha_i)}_{b_2 ,c_2} (Y_2 ; y,z,n,m,q)
$$
where 
$$
\tau = \Pi(Y) - \Pi(Y_1) -\Pi(Y_2), \quad \chi = -( \vec{b}, B^{-1} \vec{c}) + ( \vec{b}_1 , B^{-1} \vec{c}_1) +  ( \vec{b}_2 , B^{-1} \vec{c}_2) \in \rational
$$
for any choice of chamber $\alpha_i , i=\pm$.
\end{Thm}

\begin{Thm}
Let $Y_K$ be the complement of a knot $K$ in the 3-sphere $S^3$ and let $Y_{p/r}$ be a result of Dehn surgery along $K$ with slope $p/r \in \rational^{\ast}$. Assume that $Y_K$ and $Y_{p/r}$ are represented by negative definite plumbings. Then the invariants of $Y_{p/r}$ are given by
$$
\hat{Z}_{b,c}[Y_{p/r};q] = (-1)^{\tau} \mathcal{L}^{(\alpha_i ;\, p/r )}_{b,c} \lsb F^{(\alpha_i)}_K (y,z,q) \rsb,
$$
where the Laplace transform for $\alpha_+ $ chamber is
$$
\mathcal{L}^{(\alpha_+ ;\, p/r )}_{b,c} : y^{\alpha}z^{\beta}q^{\gamma} \mapsto q^{\gamma}
\begin{cases}
\sum\limits_{r_s = r_{s,min}}^{\infty}  q^{\frac{\beta(r \alpha + \epsilon r_s )}{p}}, & \text{if}\quad r \alpha + \epsilon r_s + b \in p\intg,\, r\beta + c \in p\intg\\
\sum\limits_{w_s = w_{s,min} }^{\infty} q^{\frac{\alpha(r \beta - \epsilon w_s )}{p}}, & \text{if}\quad   r \beta - \epsilon w_s + c \in p\intg ,\, r\alpha + b \in p\intg\\
0, & \text{otherwise}
\end{cases}
$$
and the Laplace transform for $\alpha_- $ chamber is 
$$
\mathcal{L}^{(\alpha_- ;\, p/r )}_{b,c} : y^{\alpha}z^{\beta}q^{\gamma} \mapsto -q^{\gamma}
\begin{cases}
\sum\limits_{w^{\prime}_s =  w^{\prime}_{s,min}  }^{\infty} q^{\frac{\beta(r \alpha - \epsilon w^{\prime}_s )}{p}}, & \text{if}\quad   r \alpha - \epsilon w^{\prime}_s + b \in p\intg ,\, r\beta + c \in p\intg\\
\sum\limits_{r^{\prime}_s =  r^{\prime}_{s,min} }^{\infty}  q^{\frac{\alpha(r \beta + \epsilon r^{\prime}_s )}{p}}, & \text{if}\quad r \beta + \epsilon r^{\prime}_s + c \in p\intg,\, r\alpha + b \in p\intg\\
0, & \text{otherwise}
\end{cases}
$$
where $r_{s,min},\, r^{\prime}_{s,min} \geq 1, \, w_{s,min},\, w^{\prime}_{s,min} \geq 0$ and $\epsilon = sign(p) (-1)^{\pi + 1}$.
\end{Thm}

\begin{Prop} Let $v$ be the number of vertices of plumbing graphs of $T(2, 2n+1)$ and $T(3, 3n+ w), w=1,2$ and $\alpha_{+}= (\alpha_1, \alpha_2 , \alpha_{v-1})$ and $\alpha_{-}$ be the good chambers for torus knots , where $\alpha_1$ corresponds to degree three vertex and the other two are associated with degree one vertices of their plumbing graphs. Their good chambers given by
$$
\alpha_{+} =  ( 1, 1 , 1),\quad \alpha_{-} = -\alpha_{+},
$$
yield a well defined (Laurent) power series $f_{m,n}(q)$. 
\end{Prop}

\begin{Conj}
Proposition 1.3 holds for all torus knots $T(s,t) \subset S^3\, (gcd(s,t)=1)$.
\end{Conj}

\begin{Conj}
Let $K \subset S^3$ be a knot and $S^{3}_{p/r}(K)$ be the result of Dehn surgery on $K$. For any choice of good chamber $\alpha_i$,
$$
\hat{Z}_{b,c}[S^{3}_{p/r}(K) ; q] = (-1)^{\tau} \mathcal{L}^{(\alpha_i ;\, p/r )}_{b,c} \lsb F^{(\alpha_i)}_K (y,z,q) \rsb,
$$
provided that the right hand side is well defined.
\end{Conj}
\noindent\textbf{Organization of the paper.} In Section 2 we review the super $\hat{Z}_{b,c}$ for closed 3-manifolds.\\
\indent In Section 3 we describe plumbed 3-manifolds and prove that the super $\hat{Z}_{b,c}$ is a topological invariant. Moreover, we describe relative $Spin^c$ structures on knot complements.\\
\indent In Section 4 we define the super $\hat{Z}_{b,c}(Y_K ; y,z,n,m,q;\alpha_i)$ for plumbed knot complements and in particular $F_K(y,z,q)$.\\
\indent In Section 5 we find chambers for torus knots and apply them to compute examples of the super $F_K(y,z,q)$.\\
\indent In Section 6 we derive the surgery formula for $F_K(y,z,q)$.\\
\indent Finally, in Section 7, we list open problems for future directions.
\newline

\noindent\textbf{Acknowledgment.} I would like to thank Heather Lee and Daren Chen for helpful explanations and Paul Orland for usage of his computer. I am grateful to Pavel Putrov for valuable comments on a draft of this paper. I wish to thank the referee for useful suggestions.

\section{Background}

We review the $q$ power series invariant of closed oriented 3-manifolds associated with the Lie superalgebra $sl(2|1)$ introduced in \cite{FP}. The physical aspects of the invariant are summarized in Appendix C.
\newline

A non-semisimple quantum invariant for a closed oriented 3-manifold $Y$ associated with $U^{(H)}_q (sl(2|1))$ at a root of unity of odd order was constructed in \cite{H}. The core ingredients of the construction are a non semi-simple ribbon category of simple finite dimensional representations of $U^{(H)}_q (sl(2|1))$ and the modified quantum dimension. The data for the quantum invariant of $Y$ consist of a root of unity of odd order $q = e^{i4\pi /l }$ ($l\geq 3$ and odd) and the degree of the Kirby color arising from a 1-cocycle, 
$$
\omega \in H^{1} (Y ; \complex / \intg \times \complex / \intg) \backslash \bigcup_{i=1}^3 H^{1} (Y ; C_i),
$$
\begin{align*}
C_1 & = \lac (X,Y) \in \complex / \intg \times \complex / \intg | 2X = 0\, \text{mod}\, 1 \rac\\
C_2 & = \lac (X,Y) \in \complex / \intg \times \complex / \intg | 2Y = 0\, \text{mod}\, 1 \rac\\
C_3 & = \lac (X,Y) \in \complex / \intg \times \complex / \intg | 2(X+Y) = 0\, \text{mod}\, 1 \rac.
\end{align*}
The non-semisimple quantum invariant is denoted by
\begin{equation}
N_{l} (Y, \omega ) \in \complex .
\end{equation}
\noindent In case of a particular class of 3-manifolds called plumbed manifolds $Y = Y(\Gamma)$~\footnote{A review of this class of manifolds is given in Seciton 3.}, it was shown in \cite{FP} that (1) decomposes into $q$-power series:
\begin{equation}
\hat{Z}^{sl(2|1)}_{b,c}[Y; q]  \in \rational + q^{\Delta_{b,c}}\intg[[q]],\qquad |q|<1,
\end{equation}
$$
(b,c) \in H_1 (Y ; \intg) \times H_1 (Y;\intg) \cong Spin^c (Y) \times Spin^c (Y),
$$
where $\Delta_{b,c} \in \rational$ and $Spin^c (Y)$ is $Spin^c$ structures on $Y$~\footnote{Its definition is a lift of the structure group $SO(3)$ of the tangent bundle $TY$ of $Y$ to $Spin^c (Y)$ group.}. This $q$ series is an analytic continuation of (1) into the complex unit disk. The decomposition of (1) is given by
\begin{multline}
N_{l} (Y(\Gamma), \omega ) = \frac{\prod\limits_{i\in V} \lb e^{i2\pi \mu^{i}_1} - e^{-i2\pi \mu^{i}_1} \rb^{\text{deg(i)-2}}}{l |Det B|}\times\\\\
\times \sum_{\substack{ \beta, \gamma \in \intg^L /B\intg^L \\ b, c \in B^{-1}\intg^L / \intg^L }} e^{i2\pi l \gamma^t B^{-1} \beta + i4\pi (b- \mu_2 )^t \gamma + i2\pi(c- (\mu_1 + \mu_2 ))^t \beta} (-1)^{\pi}\, \hat{Z}^{sl(2|1)}_{b,c}[Y(\Gamma); q]\bigg|_{q \rightarrow \zeta^2},
\end{multline}
where $\zeta = q^{1/2}$, and $(\mu^{i}_1 , \mu^{i}_2) \in \rational / \intg \times  \rational / \intg$. Furthermore, 
\begin{equation}
\hat{Z}^{sl(2|1)}_{b,c}[Y(\Gamma); q] = (-1)^{\pi}  \prod_{v \in V} \int\limits_{\Omega} \frac{dy_v}{i2\pi y_v}\frac{dz_v}{i2\pi z_v} \lb \frac{y_v - z_v}{(1-y_v)(1-z_v)} \rb^{2- \text{deg}(v_s)} \bigg|_{\alpha_i} \Theta_{b,c}(\vec{y},\vec{z},q),
\end{equation}
\begin{equation*}
\Theta_{b,c}  = \sum_{\substack{ \vec{l_1} \in B \intg^s + \vec{b} \\ \vec{l_2} \in B \intg^s + \vec{c}}} q^{(\vec{l_1}, B^{-1} \vec{l_2}) } \prod_{v \in V} y_{v}^{l_{1,v}}z_{v}^{l_{2,v}}, 
\end{equation*}
where $V$ is the vertex set of $\Gamma$, $\pi$ is the number of positive eigenvalues of $B$ and $\alpha_i$ indicates a choice of a chamber. And $\Omega$ is an integration contour tied to $\alpha_i$ chamber.\\
In contrast to $\hat{Z}_b $ associated with a Lie algebra~\cite{GPPV, P1}, the super $\hat{Z}$ (4) carries two labels $(b,c)$.
\begin{Rmk} The above integrations are equivalent to picking constant terms in the variables.
\end{Rmk}
\noindent \underline{Generic plumbing graphs} A notion of genericity of plumbing graphs was introduced in \cite{FP}. The definition states that, for a plumbing graph containing at least one vertex whose degree is greater than two, the graph does not admit a splitting $V|_{\text{deg}\neq 2} = U \sqcup W$ such that if $i\in U$ and $j\in W$, then $B^{-1}_{ij}=0$, where $V|_{\text{deg}\neq 2}$ denotes the set of vertices whose degrees are not equal to two.
\newline

\noindent \underline{Good Chambers} The integration contour $\Omega$ in (4) is equivalent to the choice of an expansion chamber $\alpha_i$. For (4) to yield a well defined power series, a (generic) plumbing graph containing at least one vertex of degree greater than two must have good chambers. The condition for the existence of good chambers is given in \cite{FP}: If there exists a vector
$$
\alpha_i = \pm 1, \qquad i \in V|_{\text{deg}\neq 2}
$$
such that 
\begin{equation}
X_{ij} : = -B^{-1}_{ij} \alpha_i \alpha_j,\quad i,j \in  V|_{\text{deg}> 2}
\end{equation}
is \textit{copositive} and
\begin{align}
B^{-1}_{ij} \alpha_i \alpha_j & \leq 0 ,\qquad \forall i \in  V|_{\text{deg}=1},\qquad j \in  V|_{\text{deg}\neq 2}\\
B^{-1}_{ij} \alpha_i \alpha_j & < 0 ,\qquad \forall i,j \in  V|_{\text{deg}=1},\qquad i \neq j
\end{align}
The matrix $X$ is called \textit{copositive} if for any vector $v$ such that $v_i \geq 0, \forall i$, with at least one $v_i \neq 0$ and have $\sum_{i,j} X_{ij} v_i v_j >0$. 
\newline

If a good chamber $\alpha$ exists for a generic plumbing graph, then there are exactly two such chambers and the domains of $y_i$ and $z_i$ corresponding to a vertex $v_i$ are given by
$$
deg(i) = 1 : \begin{cases}
|y_i|^{\alpha_i} < 1\\
|z_i|^{\alpha_i} > 1
\end{cases} \qquad
deg(i) > 2 : \bigg|\frac{y_i}{z_i}\bigg|^{\alpha_i} < 1.
$$
This translates to the following allowed expansions. For vertices $i \in V$ of degree $\text{deg}(i) = 2+ K > 2$,  expansions are
\begin{equation}
\lb \frac{(1-y_i)(1-z_i)}{y_i - z_i} \rb^K = 
\begin{cases}
(z_i -1)^K (1- y_{K}^{-1})^K \sum\limits_{r = 0}^{\infty}  \frac{(r +1)(r +2)\cdots (r +K-1)}{(K-1)!} \lb\frac{z_i}{y_i}\rb^r, & |y_i | > |z_i|\\
(1- z_{K}^{-1})^K (1- y_{K})^K \sum\limits_{r = 0}^{\infty}  \frac{(r +1)(r +2)\cdots (r +K-1)}{(K-1)!} \lb\frac{y_i}{ z_i }\rb^r, & |z_i | > |y_i|.\\
\end{cases}
\end{equation}
For vertices $i \in V$ of degree $\text{deg}(i)=1$, expansions are
\begin{equation}
\frac{y_i - z_i}{(1-y_i)(1-z_i)} = 
\begin{cases}
1+ \sum\limits_{r=1}^{\infty} y_{i}^{r} + \sum\limits_{r=1}^{\infty} z_{i}^{-r} , &  |y_i | < 1 , |z_i|>1 \\
-1- \sum\limits_{r=1}^{\infty} y_{i}^{-r} - \sum\limits_{r=1}^{\infty} z_{i}^{r} , &  |y_i| > 1 , |z_i|<1 . \\
\end{cases}
\end{equation}
Several remarks are in order.
\begin{Rmk} The other expansion domains are $|y_i|,|z_i| >1$ and $|y_i|,|z_i| < 1$. However, these are ruled out by the genericity of a plumbing graph~\cite{FP}.
\end {Rmk}
\begin{Rmk} In (2), the $\rational$ constant arises from regularizing a divergent constant. We will see the origin of this diverging constant in Sections 5 and 6.
\end {Rmk}
\begin{Rmk} The decomposition (3) was conjectured for any closed oriented 3-manifolds with $b_1 =0$ in \cite{FP}.
\end {Rmk}

\section{Plumbed manifolds}

\begin{figure}[t!]
\begin{tikzpicture}[scale=1]
\tikzstyle{every node}=[draw,shape=circle]

\draw (0,0) node[circle,fill,inner sep=1pt,label=above:$k_1 \pm 1$](){};
\draw (0,0)  -- (2,0) node[circle,fill,inner sep=1pt,label=above:$\pm 1$](){};

\draw (0,0)  -- (-1,1); 
\draw (0,0)  -- (-1,-1); 

\draw (-0.5,0.2) node[circle,fill,inner sep=1pt](){};
\draw (-0.5,0) node[circle,fill,inner sep=1pt](){};
\draw (-0.5,-0.2) node[circle,fill,inner sep=1pt](){};

\end{tikzpicture}
\quad
\begin{tikzpicture}[scale=0.9]
\tikzstyle{every node}=[draw,shape=circle]

\draw (0,0) node[circle,fill,inner sep=1pt,label=above:$ \pm 1$](){};
\draw (0,0)  -- (-2,0) node[circle,fill,inner sep=1pt,label=above:$k_1 \pm 1$](){};
\draw (0,0)  -- (2,0) node[circle,fill,inner sep=1pt,label=above:$k_2 \pm 1$](){};

\draw (2,0)  -- (3,1); 
\draw (2,0)  -- (3,-1);

\draw (-2,0)  -- (-3,1); 
\draw (-2,0)  -- (-3,-1);

\draw (-2.5,0.2) node[circle,fill,inner sep=1pt](){};
\draw (-2.5,0) node[circle,fill,inner sep=1pt](){};
\draw (-2.5,-0.2) node[circle,fill,inner sep=1pt](){};

\draw (2.5,0.2) node[circle,fill,inner sep=1pt](){};
\draw (2.5,0) node[circle,fill,inner sep=1pt](){};
\draw (2.5,-0.2) node[circle,fill,inner sep=1pt](){};

\end{tikzpicture}
\quad
\begin{tikzpicture}[scale=0.9]
\tikzstyle{every node}=[draw,shape=circle]

\draw (0,0) node[circle,fill,inner sep=1pt,label=above:$ 0$](){};
\draw (0,0)  -- (-2,0) node[circle,fill,inner sep=1pt,label=above:$k_1$](){};
\draw (0,0)  -- (2,0) node[circle,fill,inner sep=1pt,label=above:$k_2 $](){};

\draw (2,0)  -- (3,1); 
\draw (2,0)  -- (3,-1);

\draw (-2,0)  -- (-3,1); 
\draw (-2,0)  -- (-3,-1);

\draw (-2.5,0.2) node[circle,fill,inner sep=1pt](){};
\draw (-2.5,0) node[circle,fill,inner sep=1pt](){};
\draw (-2.5,-0.2) node[circle,fill,inner sep=1pt](){};

\draw (2.5,0.2) node[circle,fill,inner sep=1pt](){};
\draw (2.5,0) node[circle,fill,inner sep=1pt](){};
\draw (2.5,-0.2) node[circle,fill,inner sep=1pt](){};

\end{tikzpicture}

\begin{tikzpicture}[scale=1]
\hspace{1cm}
\rotsimeq

\hspace{5.5cm}

\rotsimeq

\hspace{5.7cm}

\rotsimeq
\end{tikzpicture}

\begin{tikzpicture}[scale=1]
\tikzstyle{every node}=[draw,shape=circle]

\draw (0,0) node[circle,fill,inner sep=1pt,label=above:$k_1$](){};

\draw (0,0)  -- (-1,1); 
\draw (0,0)  -- (-1,-1); 

\draw (-0.5,0.2) node[circle,fill,inner sep=1pt](){};
\draw (-0.5,0) node[circle,fill,inner sep=1pt](){};
\draw (-0.5,-0.2) node[circle,fill,inner sep=1pt](){};

\end{tikzpicture}
\hspace{2.3cm}
\begin{tikzpicture}[scale=0.9]
\tikzstyle{every node}=[draw,shape=circle]

\draw (0,0)  -- (-2,0) node[circle,fill,inner sep=1pt,label=above:$k_1$](){};
\draw (0,0)  -- (2,0) node[circle,fill,inner sep=1pt,label=above:$k_2 $](){};

\draw (2,0)  -- (3,1); 
\draw (2,0)  -- (3,-1);

\draw (-2,0)  -- (-3,1); 
\draw (-2,0)  -- (-3,-1);

\draw (-2.5,0.2) node[circle,fill,inner sep=1pt](){};
\draw (-2.5,0) node[circle,fill,inner sep=1pt](){};
\draw (-2.5,-0.2) node[circle,fill,inner sep=1pt](){};

\draw (2.5,0.2) node[circle,fill,inner sep=1pt](){};
\draw (2.5,0) node[circle,fill,inner sep=1pt](){};
\draw (2.5,-0.2) node[circle,fill,inner sep=1pt](){};

\end{tikzpicture}
\hspace{1.2cm}
\begin{tikzpicture}[scale=0.9]
\tikzstyle{every node}=[draw,shape=circle]

\draw (0,0) node[circle,fill,inner sep=1pt,label=above:$ k_1 + k_2 $](){};

\draw (0,0)  -- (2,1);
\draw (0,0)  -- (2,-1);

\draw (0,0)  -- (-2,1);
\draw (0,0)  -- (-2,-1);

\draw (1,0.2) node[circle,fill,inner sep=1pt](){};
\draw (1,0) node[circle,fill,inner sep=1pt](){};
\draw (1,-0.2) node[circle,fill,inner sep=1pt](){};

\draw (-1,0.2) node[circle,fill,inner sep=1pt](){};
\draw (-1,0) node[circle,fill,inner sep=1pt](){};
\draw (-1,-0.2) node[circle,fill,inner sep=1pt](){};

\end{tikzpicture}

\caption{Kirby-Neumann moves on plumbing trees. Move 1: blow up/down (left), move 2: absorption/desorption  (middle), move 3: fusion/fission (right).}
\end{figure}
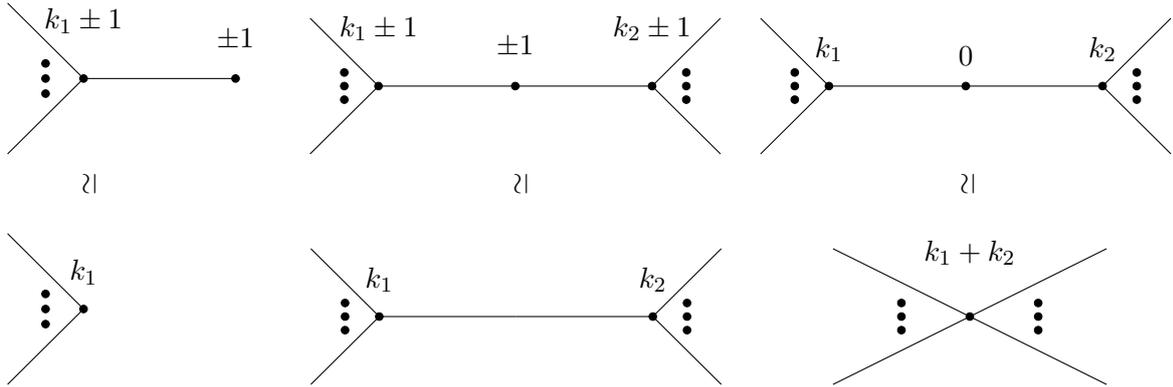

\subsection{Plumbed knot complements} We begin with a closed manifold and then move on to a knot complement. A closed oriented plumbed 3-manifold $Y$ is described by a weighted graph $\Gamma$. It consists of vertices $\lac v_i \rac$ and edges. The former carry integer weights $\lac k_i \rac$ whereas the latter carry weight $1$. This plumbing graph data is summarized by an adjacency matrix $B$, which is a symmetric and whose size is set by the number of vertices $s$ of $\Gamma$:
$$
B_{i,j} = \begin{cases}
k_i , \quad v_i = v_j \\
1 , \quad v_i , v_j \quad \text{connected}\\
0 , \quad \text{otherwise}
\end{cases}
$$
In this paper, we assume that plumbing graphs are trees. An interpretation of $\Gamma$ is that each vertex $v_i$ represents an $S^1$-bundle over $S^2$ whose Euler number is $k_i$. An edge between two vertices represents gluing two $S^1$-bundles by cutting out a $D^2$ from each base space and attaching the two $T^2$ boundaries. Another useful interpretation is a surgery link $L(\Gamma)$ obtained by replacing a vertex with a $k_i$-framed unknot and an edge with a Hopf link between two unknots. Hence $L(\Gamma)$ is always a tree link. Applying Dehn surgery (see Section 6.3) on $L(\Gamma)$ results in the same $Y$. The first homology of $Y(\Gamma)$ is
\begin{equation}
H_1 (Y(\Gamma)) \cong \intg^s / B \intg^s .
\end{equation}  
If $B$ is nondegenerate, $Y$ is a rational homology sphere. When $B$ is negative definite, we refer to $Y$ as a negative definite plumbed 3-manifold.
\newline

\indent A plumbed 3-manifold can be represented by different plumbing graphs that are related by a set of Kirby-Neumann moves in Figure 1. In \cite{Ki,N,FR}, it was shown that two plumbing graphs $\Gamma$ and $\Gamma^{\prime}$ represent the same 3-manifold $Y(\Gamma) \simeq Y(\Gamma^{\prime})$ if and only if they are related by a sequence of these moves.   
\newline

\indent A well known class of plumbed 3-manifolds consists of Seifert fibered manifolds. Their graphs are star-shaped; they consist of one central vertex of degree $\geq 2$~\footnote{Degree of a vertex is number of legs emanating from it. Degree two case is a Lens space (a special Seifert fibered manifold).} and finite number of legs attached to the central vertex. The degrees of the vertices on the legs are one or two. These legs correspond to the singular fibers of the manifold. The graph data can be summarized as follows:
$$
M \lb b \bigg| \frac{a_1}{b_1} , ..., \frac{a_n}{b_n} \rb,\qquad gcd(a_i,b_i)=1
$$
$$
e = b + \sum_{i=1}^{n} \frac{a_i}{b_i} \in \rational, \qquad ( b \in \intg )
$$
where $e$ is the Euler number, $b$ is the weight of the central vertex, $n$ is the number of singular fibers and $(a_i , b_i )$ are called the Seifert invariants. Their continued fraction expansions yield the weights of the vertices on the legs.
$$
\frac{b_i}{a_i} = 
k^{i}_1 -\cfrac{1}{k^{i}_2 -\cfrac{1}{\ddots - \frac{1}{k^{i}_s}}}.
$$
where $s$ depends on the singular fibers. A vertex attached to the central vertex has weight $-k^{i}_1$ and the last vertex on the same leg has weight $-k^{i}_s$.  
\newline
\indent For negative definite plumbed 3-manifolds, $b<0$ and $0< a_i < b_i $. It was shown in \cite{NR} that the sign of $e$ determines the positive or negative definiteness of the manifold (converse also holds).
\newline
If (7) is trivial ($H_1 =0$), $Y(\Gamma)$ is an integral homology 3-sphere. In terms of Seifert data, the $\intg HS^3$ condition is
$$
e \prod_{i=1}^{n} b_i = \pm 1.
$$
This subclass of manifolds are denoted by $\Sigma ( b_1, \cdots, b_n)$. Examples are shown in Figure 6 and 7 in Section 6.4.
\newline

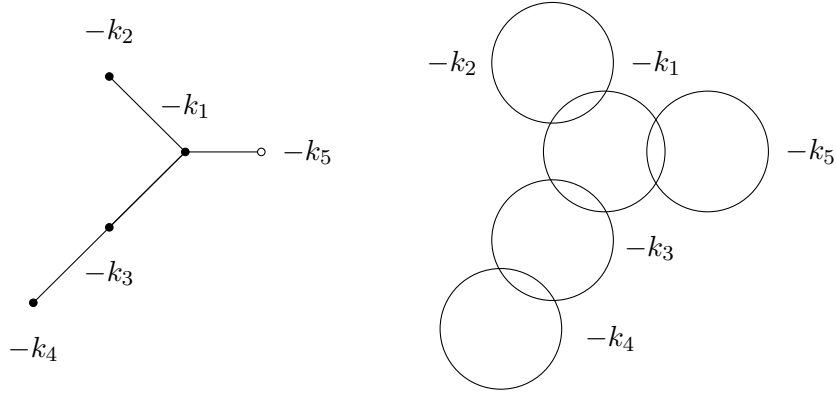
\begin{figure}[t!]
\begin{center}
\begin{tikzpicture}[scale=1]
\centering
\tikzstyle{every node}=[draw,shape=circle]

\draw (0,0) node[circle,fill,inner sep=1pt,label=above:$-k_1$](){};
\draw (0,0)  -- (1,0) node[circle,fill=white,inner sep=1pt,label=right:$-k_5$](){};
\draw (0,0)  -- (-1,1) node[circle,fill,inner sep=1pt,label=above:$-k_2$](){};
\draw (0,0)  -- (-1,-1) node[circle,fill,inner sep=1pt,label=below:$-k_3$](){};
\draw (0,0)  -- (-2,-2) node[circle,fill,inner sep=1pt,label=below:$-k_4$](){};

\end{tikzpicture}
\qquad
\begin{tikzpicture}[scale=0.8]
\centering

\draw (0,0) circle (1cm);
\draw (-0.85,1.47) circle (1cm);
\draw (-0.85,-1.47) circle (1cm);
\draw (-1.7,-2.94) circle (1cm);
\draw (1.7,0) circle (1cm);

\node at (0.85,1.47) {$-k_1$};
\node at (-2.5,1.47) {$-k_2$};
\node at (0.75,-1.6) {$-k_3$};
\node at (0.1,-3.1) {$-k_4$};
\node at (3.4,0) {$-k_5$};

\end{tikzpicture}

\end{center}
\caption{A plumbing graph $\Gamma$ of a knot $\subset S^3$ (left) and corresponding surgery link $L(\Gamma)$. The linking between two link components is the Hopf link. This link diagram can be transformed into a knot diagram through the Kirby moves.}
\end{figure}

Plumbed knot complements, more generally, plumbed 3-manifolds with a torus boundary, are represented by a weighted graph $\Gamma$ with one distinguished vertex $v_{\ast}$~\cite{GM}. This vertex represents the torus boundary. We are interested in the case where the degree of $v_{\ast}$ is one. From the viewpoint of the surgery link $L(\Gamma)$ described above, an unknot corresponding to $v_{\ast}$ acts as a spectator during the surgery operation. Furthermore, removing $v_{\ast}$ and the edge connecting it to $\Gamma$ represents an ambient plumbed 3-manifold $Y(\hat{\Gamma})$. 
\newline
\indent Additional data describing a knot is the framing, which takes values in $\intg$. Roughly speaking, this value characterizes the twisting of a longitude of the knot around the knot. This information is captured by the weight $k_{v_{\ast}}$ of $v_{\ast}$. This is called the \textit{graph framing}. Therefore, the complement of a plumbed knot in $Y(\hat{\Gamma})$ is specified by $(\Gamma , v_{\ast})$. A simple example is shown in Figure 2. The Neumann moves in Figure 1 apply to plumbing graphs of knots with the condition that the vertices of the graph need to be regular vertices. Throughout this paper, we will focus on plumbing graphs whose $B$ is (weakly) negative definite.
\begin{Defn} (Gukov-Manolescu [\cite{GM}]) Let $Y_K = Y_K (\Gamma,v_{\ast})$ be a plumbing tree consisting of $s$ number of vertieces and $v_{\ast}$ be a distinguished vertex. The pair $ (\Gamma,v_{\ast})$ is called called weakly negative definite if if the corresponding matrix $B$ is invertible, and $B^{-1}$ is negative definite on the subspace of $\intg^s$ spanned by the non-distinguished vertices of degree $\geq$ 3.  
\end{Defn}

\subsection{Invariance}

We show that the super $\hat{Z}_{b,c}$ is a topological invariant of plumbed 3-manifolds. We begin with the existence of good chambers.

\begin{Lem} The existence of good chambers is preserved under the Kirby-Neumann moves shown in Figure 1.
\end{Lem}
\noindent \textit{Proof}. For move 1, we begin with the top graph $\Gamma$ consisting of $s$ be the number of vertices and its adjacency matrix $B$ admitting good chambers. Let $v_s$ be the degree one vertex with framing $\pm 1$ and $B^{\prime}$ be an adjacency matrix of the bottom graph. In case of degree of the vertex $v_{s-1}$ of $\Gamma$ is greater than two, after blow down, $X^{\prime}_{MN}$ of the bottom graph is copositive because the submatrices of $B$ and of $B^{\prime}$ corresponding to the left of $v_{s-1}$ are the same. Furthermore, $v_s$ is only connected to $v_{s-1}$. Hence it does not affect the part of the top graph left of $v_{s-1}$.
\newline

\noindent In case of degree of the vertex $v_{s-1}$ of $\Gamma$ is two and degree of $v_{s-1}$ is greater than two, there are two subcases. If $B^{\prime\,-1}_{s-2,s-1}=0$, then (6) is fulfilled. If $B^{\prime\,-1}_{s-2,s-1}\neq 0$, then $\alpha^{\prime}_{s-2} = \alpha_{s-2}$ and $\alpha^{\prime}_{s-1}$ is determined by the sign of $B^{\prime\,-1}_{s-2,s-1}$ in (6). Moreover, (5) and (7) for other vertices are not affected by $\alpha^{\prime}_{s-1}$ since $v_{s-1}$ is only connected to $v_{s-2}$.
Copositivity of $X^{\prime}_{MN}$ is ensured by the same reason as the above case.  
\newline

\noindent For move 2, degree of the middle vertex with framing $\pm 1$ in the top graph is two, therefore it does not influence (5-7).
\newline

\noindent For move 3, We denote the middle vertex by $v_0$, the left vertex by $v_1$ and the right vertex by $v_2$ of the top graph $\Gamma$. If degrees of $v_1$ and $v_2$ are two, then the fusion of $v_1$ and $v_2$ does not affect (5-7) because the vertex with framing $k_1 + k_2$ has degree two. 	If degrees of $v_1$ and $v_2$ are greater two, then the fusion of $v_1$ and $v_2$ transfers the copositive property of $X_{IJ}$ to $X^{\prime}_{MN}$ of the bottom graph since the submatrices of $B$ and $B^{\prime}$ corresponding to the left and right parts are the same. And $v_0$ is disconnected from the left and right sides of $v_1$ and $v_2$, respectively. If degree of $v_1$ is two and degree of $v_2$ are greater two, after fusion, copositive property of $X_{IJ}$ is inherited to $X^{\prime}_{MN}$ of the bottom graph since degree two vertices are excluded in (5) and hence does not affect the left side of $v_1$ and the right side of $v_2$. Moreover, degree of $v_2$ remains the same.

\begin{Prop} The super $q$-series $\hat{Z}_{b,c}$ defined in (4) is invariant under the Kirby-Neumann moves shown in Figure 1.
\end{Prop}
\noindent \textit{Proof}. Consider the move 2 with the $-1$ signs in the graphs. Let 
\begin{equation}
\vec{m} = ( \vec{m}_L , \vec{m}_R )\qquad \vec{n} = ( \vec{n}_L , \vec{n}_R ) \in \intg^{s},
\end{equation}
be the lattice vectors of bottom graph, where $s$ is the number of vertices in the bottom graph and $\vec{m}_L ,\, \vec{m}_R$ are left and right sides of the graph, respectively. We denote its adjacency matrix by $B^{\prime}$ and the variables of the vertex having the sign with $(z_0 , y_0 )$. In the theta function of $B^{\prime}$, there is an extra factor $z^{m^{\prime}_0}_0, y^{n^{\prime}_0}_0 $. Recall that the integrations in (4) pick out constant terms. Hence only $m^{\prime}_0 = n^{\prime}_0 = 0$ contributes. Then from $\vec{m}$ and $\vec{n}$
we can obtain lattice vectors for the top graph.
$$
\vec{m}^{\prime} = ( \vec{m}_L ,0, \vec{m}_R )\qquad \vec{n}^{\prime} = ( \vec{n}_L ,0, \vec{n}_R ) \in \intg^{s+1}.
$$
From linear algebra, the exponents of $q$ are the same
$$
(\vec{m}^{\prime} , B^{\prime\, -1} \vec{n}^{\prime}) = (\vec{m} , B^{-1} \vec{n}).
$$
\newline
\indent In case of $+1$ signs, $\pi^{\prime} = \pi + 1 $. Given (11), the lattice vectors for the top graph are
$$
\vec{m}^{\prime} = ( \vec{m}_L ,0, -\vec{m}_R )\qquad \vec{n}^{\prime} = ( \vec{n}_L ,0, -\vec{n}_R ).
$$
have the same $q$ exponent as the bottom graph. Because of the above sign changes, the variables in the lattice theta function of the vertices of the right side of the top graph need to be change to $z_v \rightarrow z_{v}^{-1} , \,  y_v \rightarrow y_{v}^{-1} $. This results in an extra $(-1)^r$ factor, where
$$
r= \sum_{v \in \Gamma^{\prime}_R} 2- \text{deg}(v).
$$
This value is odd. Thus this $-1$ sign compensates the additional sign from $\pi^{\prime}$.
\newline

Next, we consider the move 1 with $-1$ sign. Let the lattice vectors for the bottom graph be
$$
\vec{m} = ( \vec{m}_L , m_1 )\qquad \vec{n} = ( \vec{n}_L , n_1  ),
$$  
where $m_1$ and $n_1$ are in the entry for the vertex with weight $k_1$. Corresponding lattice vectors for the top graph are
$$
\vec{m}^{\prime} = ( \vec{m}_L , m_1 - m_0 , m_0 )\qquad \vec{n}^{\prime} = ( \vec{n}_L ,  n_1 -  n_0 , n_0 ).
$$
The super $\hat{Z}$ of the top graph has extra factors
$$
\lb \frac{y_1 - z_1}{(1- y_1)(1-z_1)} \rb^{-1} \frac{y_0 - z_0}{(1- y_0)(1-z_0)} z_{1}^{- m_0} y_{1}^{- n_0} z_{0}^{ m_0} y_{1}^{ n_0} 
$$
We expand the edge term for $(y_0,z_0)$ in one of the chambers (9). The contributing values of $( m_0, n_0 )$ are
\begin{equation}
( m_0, n_0 ) \in \lac (0,0) , (0, -r_0 ), (r_0 , 0) \, |\, r_0 \in \intg_+ \rac.
\end{equation}
This implies that $z_{1}^{- m_0} y_{1}^{- n_0}$ term yields 
$$
1+ \sum_{r_0 = 1}^{\infty} z_{1}^{-r_0} + \sum_{r_0 = 1}^{\infty} y_{1}^{r_0} = \frac{y_1 - z_1}{(1- y_1)(1-z_1)},
$$
where (9) is used. This cancels the above edge term for $(y_1,z_1)$ . Similar cancellation occurs for the other chamber in (9). We next compare the exponents of $q$,
$$
\lb \vec{m}^{\prime} , B^{\prime\,-1} \vec{n}^{\prime} \rb - (\vec{m} , B^{-1} \vec{n}) = - m_0 n_0
$$
Observe that all elements of (12) has $m_0$ or $n_0$ is zero. Thus the powers of $q$ of the top and the bottom graphs are match.\\
\indent In case of move 1 with $+1$ sign, we use
$$
\vec{m}^{\prime} = ( \vec{m}_L , m_1 + m_0 , m_0 )\qquad \vec{n}^{\prime} = ( \vec{n}_L ,  n_1 +  n_0 , n_0 ).
$$
\newline
\indent For move 3, the top graph has $\pi^{\prime} = \pi + 1$ due to an extra positive eigenvalue. We denote the middle vertex by $v_0$, the left vertex by $v_1$ and the right vertex by $v_2$. The super $\hat{Z}$ of the graph contain the term
\begin{equation}
\lb \frac{y_1 - z_1}{(1-y_1)(1-z_1)} \rb^{2-  \text{deg} (v_1)}  \lb \frac{y_2 - z_2}{(1-y_2)(1-z_2)} \rb^{2-  \text{deg} (v_2)}  z_{0}^{m_0} y_{0}^{n_0} z_{1}^{m_1} y_{1}^{n_1} z_{2}^{m_2} y_{2}^{n_2}.
\end{equation}
The integrations over $z_{0}$ and $y_{0}$ imply that $m_0 = n_0 =0$. Let the lattice vectors be
$$
\vec{m}^{\prime} = ( \vec{m}_L , m_1 , 0, m_2 , \vec{m}_R  ), \qquad \vec{n}^{\prime} = ( \vec{n}_L , n_1 , 0, n_2 , \vec{n}_R  ) \in \intg^{s+2},
$$ 
where $\vec{m}_L$ and $\vec{m}_R $ are associated with left and right part of the graph excluding $v_1$ and $v_2$. Similarly for $\vec{n}_L$ and $\vec{n}_R$. Matching of the $q$ exponents between the top and the bottom graphs $\lb \vec{m}^{\prime} , B^{\prime\,-1} \vec{n}^{\prime} \rb = (\vec{m} , B^{-1} \vec{n})$ requires the lattice vectors of the latter graph be
$$
\vec{m} = ( \vec{m}_L , m_1 - m_2 , -\vec{m}_R  ), \qquad \vec{n} = ( \vec{n}_L , n_1 - n_2 , -\vec{n}_R  ) \in \intg^{s}.
$$
This implies that we need to invert the variables associated with the vertices in the right part of the top graph and those corresponding to $v_1$ and $v_2$.
$$
z_v \rightarrow z_{v}^{-1}, \quad  y_v \rightarrow y_{v}^{-1} ,\quad z_2 \rightarrow z_{1}^{-1}, \quad  y_2 \rightarrow y_{1}^{-1}.
$$
After defining $z_r : = z_1, y_r := y_1$ for the central vertex $v_r$ of the bottom graph, (13) becomes
$$
\pm \lb \frac{y_r - z_r}{(1-y_r)(1-z_r)} \rb^{2-  \text{deg} (v_r)} z_{r}^{m_1 - m_2 } y_{r}^{n_1 - n_2 },
$$
where 
$$
\lb 2 - \text{deg}(v_1) \rb + \lb 2 - \text{deg}(v_2) \rb  = 2 -  \text{deg} (v_r)
$$
is used. The minus sign corresponds to the case of even number of degree one vertices whereas the plus sign for the case of odd number of degree one vertices. In the latter case, an extra minus sign comes from inverting $z_v$ and $y_v$ of the degree one vertices. Thus, in both cases, there is a minus sign that cancels the minus sign from $\pi^{\prime}$. Thus we arrive at the super $\hat{Z}$ of the bottom graph.

\subsection{$Spin^c$ Structures on knot complements}

In the case of knot complements, the labels $(b,c)$ of the super $\hat{Z}$ are elements of  $H_1 (Y_K) \times H_1 (Y_K)$, which is affinely isomorphic to $Spin^c (Y_K, \ptl Y_K) \times Spin^c (Y_K, \ptl Y_K)$.  We describe relative $Spin^c$ structures on plumbed knot complements~\cite{GM}.
\newline
\indent Let $Y_K = Y_K(\Gamma, v_{\ast})$ be the complement of a knot $K$ in $S^3$ represented by $(\Gamma, v_{\ast})$. The graph consists of $s$ number of vertices, where $v_s = v_{\ast}$. We denote its adjacency matrix by $B$. Then we have
$$
H^2 (Y_K, \ptl Y_K) \cong H_1 (Y_K) \cong \intg^s / B \intg^{s-1} 
$$
where $\intg^{s-1} =  \intg^{s-1} \times \lac 0 \rac \subset  \intg^{s}$.\\
\indent We let $\vec{e}_i ,\, i=1,\cdots ,s$ be the basis vectors of $\intg^{s}$. The meridian and longitude of the $T^2$ boundary of $Y_K$ are $\vec{e}_s$ and $B\vec{e}_s$, respectively.
$$
H_1 (T^2) \cong Span\left\langle \vec{e}_s , B\vec{e}_s \right\rangle \subset  \intg^{s}.
$$
The action of $H_1 (T^2)$ on $ H_1 (Y_K) $ is given by adding multiples of $\vec{e}_s$ and $B\vec{e}_s$. The above two identifications can be combined into
$$
Span\left\langle \vec{e}_s , B\vec{e}_s \right\rangle \hookrightarrow    \intg^{s}  \twoheadrightarrow  \intg^s / B \intg^{s-1}.
$$
The relative $Spin^c$ structures on $Y_K$ are
$$
Spin^c (Y_K, \ptl Y_K) \cong  2 \intg^{s} + \vec{\delta} / (2B\intg^{s-1}).
$$
The $Spin^c$ structures on $Y_K$ are
$$
Spin^c (Y_K) \cong  2 \intg^{s} + \vec{\delta} / (Span\left\langle 2\vec{e}_s , 2B\vec{e}_s \right\rangle +  2B\intg^{s-1}).
$$

\section{A supergroup series invariant of plumbed knot complements}

Motivated by the idea of partial surgery in \cite{GM}, we define a series invariant $\hat{Z}_{b,c}$ for plumbed knot complements in this section.

\subsection{A partial surgery formula}

\indent We use the surgery link $L(\Gamma)$ interpretation of plumbing graphs of knots in $S^3$ to write down a partial surgery formula. Recall that the torus boundary of a plumbed knot complement is represented by a distinguished vertex $v_s$ in the plumbed graph, which is depicted as an open circle as in Figure 2. Such a vertex carries four variables, which we denote by
$$ y= y_s,\quad z=z_s,\quad n=n_s , \quad m=m_s. $$
We apply partial surgery on $L(\Gamma)$ by integrating over its link components, except the link component corresponding to $v_s$. 
\begin{Defn} For a plumbed knot complement $Y_K =(\Gamma, v_s)$ with a generic $\Gamma$ admitting good chambers $\alpha_i$, $i=\pm$ and a weakly negative definite $B$,we define a super series invariant in $\alpha_i$ chamber by
\end{Defn}
\begin{align}
\hat{Z}^{sl(2|1)}_{b,c}(Y_K; y , z , n,m,q; \alpha_i) & = (-1)^{\pi} \lb \frac{y-z}{(1-y)(1-z)} \rb^{1-\text{deg}(v_s)} \prod_{\substack{ v \in V \\ v \neq v_s } } \int\limits_{\Omega} \frac{dy_v}{i2\pi y_v}\frac{dz_v}{i2\pi z_v}\\
& \times \lb \frac{y_v - z_v}{(1-y_v)(1-z_v)} \rb^{2- \text{deg}(v_s)}  \Theta_{b,c}(\vec{y},\vec{z},q),\nonumber\\
\Theta_{b,c}(\vec{y},\vec{z},q) & = \sum_{\substack{ \vec{l_1} \in B \intg^s + \vec{b} \\ \vec{l_2} \in B \intg^s + \vec{c}}} q^{(\vec{l_1}, B^{-1} \vec{l_2}) } \prod_{v \in V} y_{v}^{l_{1,v}}z_{v}^{l_{2,v}},\nonumber
\end{align}
where the last components of $\vec{n} \in \vec{l_1} = B \vec{n} + \vec{b} $  and  $\vec{m} \in \vec{l_2} = B \vec{m} + \vec{c} $ are  $n$ and $m$, respectively. The good chamber expansions for $(y_v,z_v)$ are given in (8) and (9).
\begin{Rmk} In the case of a nondegenerate $B$, the existence conditions of good chambers for a plumbed knot complement are the same as (5-7), except that the distinguished vertex is excluded from them.
\end{Rmk}


\noindent The exponent of the prefactor $(y-z)/((1-y)(1-z))$ is $1-\text{deg}(v_s)$ for the purpose of gluing two knot complements (see Section 6). As in the closed oriented manifold case (4), the integration contour $\Omega$ corresponds to a choice of a good chamber $\alpha_i$. 
\newline

\indent Relative $\text{Spin}^c$ structures of plumbed knot complements carry a conjugation symmetry, as seen in (14):
\begin{equation}
\hat{Z}^{(\alpha_+)}_{-b,-c}[Y_K; y , z , n,m,q] = - \hat{Z}^{(\alpha_-)}_{b,c}[Y_K; y^{-1} , z^{-1} , -n,-m,q ].
\end{equation}
This symmetry exchanges the two chambers $\alpha_+$ and $\alpha_-$ as the domains of $y$ and $z$ are switched. This is in contrast to the case of regular Lie groups~\cite{GM}; the series invariant is invariant under the above symmetry transformation because there is no notion of expansion chambers.
\begin{Rmk} We specify chambers $\alpha_{\pm}$ as a superscript or in brackets of the super $\hat{Z}$.
\end{Rmk}

\noindent The complete series invariant is given by the sum of the two chamber contributions
\begin{equation}
\hat{Z}_{b,c} [Y_K; y , z , n,m,q] = \hat{Z}^{(\alpha_+)}_{b,c}[Y_K; y , z , n,m,q] + \hat{Z}^{(\alpha_-)}_{b,c}[Y_K; y , z , n,m,q].
\end{equation}
The relative $Spin^c$ conjugation symmetry in (15) translates into a Weyl symmetry in (16). Hence, (16) is manifestly Weyl symmetric in $y$ and $z$.
\newline

\noindent \underline{Degenerate $B$} For some knots in $S^3$, their adjacency matrices $B$'s are non-invertible. In such cases, the lattice theta function in (14) needs to be modified. We let $\vec{b}=B\vec{g}$ and $\vec{c}=B\vec{w}$ for some $\vec{g}, \vec{w} \in \intg^{s-1}$. Then we have $\vec{l}_1 = B(\vec{n} + \vec{g}),\, \vec{l}_2 = B(\vec{m} + \vec{w})$. The theta function becomes
\begin{equation}
\Theta_{b,c} = q^{(\vec{g},B\vec{w})} \sum_{\vec{n} \in \intg^s} \sum_{\vec{m} \in \intg^s} q^{(\vec{n},B\vec{m}) + (\vec{n},B\vec{w})+ (\vec{m},B\vec{g})}  \prod_{v \in V} y_{v}^{l_{1,v}(\vec{n})}z_{v}^{l_{2,v}(\vec{m})}.
\end{equation}
We will see in Section 5.2 that there exist expansion chambers for torus knots in which the exponent of $q$ in (17) is bounded below (i.e. good chambers).
\newline

\begin{Prop} The super $\hat{Z}_{b,c}$ of the plumbed knot complements (14) is invariant under the Kirby-Neumann moves in Figure 1.
\end{Prop}
\noindent \textit{Proof.} The proof is same as that of Proposition 3.3.

\subsection{The solid torus} 

\begin{figure}[t]
\begin{center}
\begin{tikzpicture}[scale=1]
\centering
\tikzstyle{every node}=[draw,shape=circle]

\draw (0,0) node[circle,fill=white,inner sep=1pt,label=above:$-k_1$](){};
\draw (0,0)  -- (1,0) node[circle,fill,inner sep=1pt,label=above:$-k_2$](){};
\draw (1,0) -- (5,0) node[circle,fill,inner sep=1pt,label=above:$-k_s$](){};

\draw (2,0.1) node[circle,fill,inner sep=1pt,label=above:$$](){};
\draw (3,0.1) node[circle,fill,inner sep=1pt,label=above:$$](){};
\draw (4,0.1) node[circle,fill,inner sep=1pt,label=above:$$](){};

\end{tikzpicture}

\end{center}
\caption{Plumbing graphs of the solid torus $S_{p/r}$. The distinguished vertex is the first vertex as indicated by an open circle. The ellipsis indicates intermediate vertices on the leg whose framing coefficients are determined by the continued fraction expansion of $p/r$ in Section 3.1.}
\end{figure}
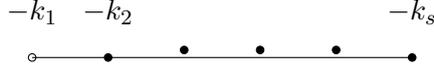
We compute the super $\hat{Z}$ for the solid torus $S_{p/r}$ for $r\neq 0$ and $gcd(p,r)=1$. Its plumbing graph is given by a linear plumbing shown in Figure 3. The simplest case is a graph having a single distinguished vertex  with $r=1$. It represents a $p$-framed unknot $U_p\, (p\neq 0)$. From (14), we get
\begin{equation}
\hat{Z}_{b,c} [U_p ; y , z , n,m,q] = \pm sign(p) \frac{y-z}{(1-y)(1-z)} q^{\frac{(pn + b)(pm + c)}{p}} y^{pn+b} z^{pm+c}.
\end{equation}
The rational function in (18) can be expanded in two ways depending on the choice of good chambers (9). By setting  $b$ and $c$ to zero in (18), we can take the $p=0$ limit. Then we obtain the 0-framed unknot result.
\begin{equation}
\hat{Z}_{0,0} [U ; y , z , q] =  \frac{y-z}{(1-y)(1-z)}.
\end{equation}

We next move on to a generic (knotted) solid torus having at least two vertices in its graph $(v \geq 2)$. According to Figure 3, there is one degree one vertex $v_s$ that contributes to $\hat{Z}$ of the solid torus. In the $\alpha_+$ chamber, we have
$$
\lb \sum_{r_s =1}^{\infty} y_s^{r_s} + \sum_{w_s = 0}^{\infty} z_s^{-w_s} \rb y_s^{l_{1,s}}z_s^{l_{2,s}}
$$
from (14). Recalling Remark 2.1, integration over the variables $y_s$ and $z_s$ implies that
\begin{align*}
l_{1,s} & = - r_s,& l_{2,s} &  =0\\
l_{2,s} & =   w_s,& l_{1,s} & =0
\end{align*}
Then $\vec{l_1}$ and $\vec{l_2}$ have the following components.
\begin{align}
\begin{split}
\Lambda_{b,c}^{-,0} & = \lac \lb \vec{l_1} = ( r_1, 0,...,0,-r_s) , \vec{l_2} = ( g_1, 0,...,0,0) \rb \bigg| r_s \in \intg_+ , r_1, g_1 \in \intg \rac\\
\Lambda_{b,c}^{0,+} & = \lac \lb \vec{l_1} = ( d_1, 0,...,0,0) , \vec{l_2} = (w_1, 0,...,0,w_s) \rb \bigg| w_s \in \intg_{\geq 0} , d_1, w_1 \in \intg \rac.
\end{split}
\end{align}
Hence, the $\hat{Z}$ in this chamber is
\begin{equation}
\hat{Z}^{(\alpha_+)}_{b,c}[S_{p/r}; y_1 , z_1 , n_1 , m_1 , q] = (-1)^{\pi} \lb \sum_{\vec{l_i} \in \Lambda_{b,c}^{-,0}} y_1^{l_{1,1}}z_1^{l_{2,1}} q^{\vec{l_1} B^{-1} \vec{l_2}} +  \sum_{\vec{l_i} \in \Lambda_{b,c}^{0,+}} y_1^{l_{1,1}}z_1^{l_{2,1}} q^{\vec{l_1} B^{-1} \vec{l_2}} \rb,
\end{equation}
where $\vec{l_1} = B \vec{n} + \vec{b}, \vec{l_2} = B \vec{m} + \vec{c} $ for some $\vec{n} = (n_1, n_2,...,n_s)$ and  $\vec{m} = (m_1, m_2,...,m_s)$.\\
\noindent The adjacency matrix of $S_{p/r}$ in Figure 3 is
$$
B= \begin{pmatrix}
k_1 & 1   & 0 & \dots  &  \dots &  \dots\\
1   & k_2 & 1 & 0 &  \dots &  \dots\\
0   & 1   & k_3 & 0 &  \dots &  \dots\\
0   & \dots &  \dots   &  \ddots &  \dots &  \dots\\
0   & \dots &  \dots   & \dots &  k_{s-1} & 1 \\
0   & \dots &  \dots & \dots &  1 & k_s \\
\end{pmatrix}
$$
Its determinant is $\pm p$. We next calculate the exponents of $q$ in (21). 
\begin{align*}
\vec{l_1} B^{-1} \vec{l_2} & = ( r_1 , \vec{0}, -r_s )^t  B^{-1} (g_1 , \vec{0} )\\
                           & =  \frac{g_1}{p} \lb r r_1  - \epsilon r_s \rb\\
\vec{l_1} B^{-1} \vec{l_2} & = ( d_1 , \vec{0} )^t  B^{-1} (w_1 , \vec{0}, w_s )\\
                           & =  \frac{d_1}{p} \lb r w_1  + \epsilon w_s \rb
\end{align*}
where $B^{-1}_{11} = r/p , B^{-1}_{s1}= \epsilon /p$ and $\epsilon = sign(p) (-1)^{\pi + 1}$ are used. After substitutions into (21), we arrive at
\begin{equation}
\hat{Z}^{(\alpha_+)}_{b,c}[S_{p/r}; y_1 , z_1 , n_1 , m_1 , q] = (-1)^{\pi} \lb \sum_{ \Lambda_{b,c}^{-,0}} y_1^{r_1} z_1^{g_1} q^{\frac{g_1}{p} \lb r r_1  - \epsilon r_s \rb} +  \sum_{\Lambda_{b,c}^{0,+}} y_1^{d_1}z_1^{w_1} q^{ \frac{d_1}{p} \lb r w_1  + \epsilon w_s \rb} \rb.
\end{equation}
We note that the each summation in the above is multiple summations whose ranges are given by (20).
\newline
\indent Applying the same method to the second chamber $\alpha_-$, we obtain
\begin{equation}
\hat{Z}^{(\alpha_-)}_{b,c}[S_{p/r}; y_1 , z_1 , n_1 , m_1 , q] = (-1)^{\pi +1} \lb \sum_{ \Lambda_{b,c}^{+,0}} y_1^{w^{\prime}_1} z_1^{h_1} q^{\frac{h_1}{p} \lb r w^{\prime}_1  + \epsilon w^{\prime}_s \rb} +  \sum_{\Lambda_{b,c}^{0,-}} y_1^{u_1}z_1^{r^{\prime}_1} q^{ \frac{u_1}{p} \lb r r^{\prime}_1  - \epsilon r^{\prime}_s \rb} \rb,
\end{equation}
where
\begin{align}
\Lambda_{b,c}^{+,0} & = \lac \lb \vec{l_1} = ( w^{\prime}_1, 0,...,0, w^{\prime}_s) , \vec{l_2} = ( h_1, \vec{0} ) \rb \bigg| w^{\prime}_s \in \intg_{\geq 0}, w^{\prime}_1, h_1 \in \intg \rac \nonumber\\
\Lambda_{b,c}^{0,-} & = \lac \lb \vec{l_1} = ( u_1, \vec{0}) , \vec{l_2} = (r^{\prime}_1, 0,...,0,-r^{\prime}_s) \rb \bigg| r^{\prime}_s \in \intg_{+} , u_1, r^{\prime}_1 \in \intg \rac
\end{align}
Each summation in the above is multiple summations whose ranges are given by (24). 
We next express (22) in terms of $m$ and $n\, (m=m_1 , n = n_1)$. The values of $r_1, r_s$ and $g_1$ in $\alpha_+$ chamber are related via $m$ and $n$ in the following way.
\begin{align}
r r_1 & = pn + b + \epsilon r_s   & r g_1 & = pm + c \nonumber\\
r w_1 & = pm + c - \epsilon w_s   & r d_1 & = pn + b
\end{align}
Substituting them into (22) yields
\begin{align*}
\hat{Z}^{(\alpha_+)}_{b,c}[S_{p/r}; y_1 , z_1 , n , m , q] & = (-1)^{\pi } q^{\frac{(pn+b)(pm+c)}{pr}} \lb \sum_{r_s =1 }^{\infty} y_1^{\frac{pn+b+\epsilon r_s}{r}} z_1^{\frac{pm+c}{r}}  + \sum_{w_s =0 }^{\infty} y_1^{\frac{pn+b}{r}} z_1^{\frac{pm+c -\epsilon w_s }{r}} \rb\\
& = (-1)^{\pi } q^{\frac{(pn+b)(pm+c)}{pr}} \lb y_1^{\frac{pn+b}{r}} z_1^{\frac{pm+c}{r}} + \sum_{j =1 }^{\infty} y_1^{\frac{pn+b+\epsilon j}{r}} z_1^{\frac{pm+c}{r}} + y_1^{\frac{pn+b}{r}} z_1^{\frac{pm+c -\epsilon j }{r}} \rb.
\end{align*}
In case of $\alpha_-$ chamber, we have
\begin{align*}
r w^{\prime}_1 & = pn + b - \epsilon w^{\prime}_s   & r h_1 & = pm + c \\
r r^{\prime}_1 & = pm + c + \epsilon r^{\prime}_s   & r u_1 & = pn + b.
\end{align*}
After substitutions	into (23), we arrive at
\begin{align*}
\hat{Z}^{(\alpha_-)}_{b,c}[S_{p/r}; y_1 , z_1 , n , m , q] & = (-1)^{\pi +1}  q^{\frac{(pn+b)(pm+c)}{pr}} \lb \sum_{ w^{\prime}_s =0 }^{\infty} y_1^{\frac{pn+b- \epsilon w^{\prime}_s}{r}} z_1^{\frac{pm+c}{r}} + \sum_{r^{\prime}_s =1 }^{\infty} y_1^{\frac{pn+b}{r}} z_1^{\frac{pm+c + \epsilon r^{\prime}_s }{r}} \rb\\
& =  (-1)^{\pi +1}  q^{\frac{(pn+b)(pm+c)}{pr}} \lb  y_1^{\frac{pn+b}{r}} z_1^{\frac{pm+c}{r}} + \sum_{ j = 1 }^{\infty} y_1^{\frac{pn+b- \epsilon j}{r}} z_1^{\frac{pm+c}{r}} + y_1^{\frac{pn+b}{r}} z_1^{\frac{pm+c + \epsilon j }{r}} \rb.
\end{align*}
Therefore,
\begin{equation}
\hat{Z}_{b,c}[S_{p/r}; y_1 , z_1 , n , m , q] = \hat{Z}^{(\alpha_+)}_{b,c}[S_{p/r}; y_1 , z_1 , n , m , q]  + \hat{Z}^{(\alpha_-)}_{b,c}[S_{p/r}; y_1 , z_1 , n , m , q].
\end{equation}
Under (15), it is straightforward to check that the two terms in the right hand side of (26) exchange.

\subsection{The boundary action}

For $\hat{Z}$ associated with a Lie algebra $sl(2)$, a useful simplification arose due to the $H_1 (T^2)$ action of the boundary torus on the label of the $\hat{Z}$, which takes value in $Spin^c (Y_K , \ptl Y_K)$ of a plumbed knot complement~\cite{GM}
\newline
\indent In case of $sl(2|1)$, we have an action of $H_1 (T^2) \times H_1 (T^2)$ on the labels of super $\hat{Z}_{b,c}$ taking values in $H_1 (Y_K) \times H_1 (Y_K)$~\footnote{This is a consequence of $H_1 (T^2) \times H_1 (T^2)$ action on the vector space assigned to $T^2$ (see Section 6.4 for details).}. This action entails the following consequences. 
\begin{Prop} Let $Y_K = Y(\Gamma, v_s)$ be a negative definite plumbed knot complement. Then for any $(b,c) \in H_1(Y_K) \times H_1(Y_K)$ and $(\gamma , \eta) \in H_1 (T^2) \times H_1 (T^2)$, we have
$$
A_{\gamma , \eta} \hat{Z}_{b,c} \cong  \hat{Z}_{b + g(\gamma),c + g (\eta)}.
$$
\end{Prop}
\noindent \textit{Proof.} Let $\vec{b}$ and $\vec{c}$ be vector representatives of $(b,c)$ in $\intg^s$. The action of the meridian components of $(\gamma , \eta)$ amounts to adding $\vec{e_s}$ to $\vec{b}$ and $\vec{c}$. This shifts 
$$
b_s \mapsto  b_s + 1 \qquad c_s \mapsto  c_s + 1
$$
To analyze effects, we substitute $\vec{l}_1$ and $\vec{l}_2$ expressions into the summand in (14). Then we find that the $q$ exponent $(\vec{n},\vec{c}) + (\vec{m},\vec{b}) + (\vec{b}, B^{-1}\vec{c})$ is shifted. Furthermore, extra multiplicative $y_s$ and $z_s$ factors appear. Overall, the above actions result in a multiplication by a $q$ monomial.
\newline
\indent The action of the longitude components of $(\gamma , \eta)$ is done by adding $B\vec{e}_s$ to $\vec{b}$ and $\vec{c}$. Consequently,
$$
\vec{l}_1 = B \vec{n} + \vec{b}  = B (\vec{n} - \vec{e}_s )  + (\vec{b} + B\vec{e}_s ),\quad \vec{l}_2 = B \vec{m} + \vec{c}  = B (\vec{m} - \vec{e}_s )  + ( \vec{c} + B\vec{e}_s )
$$
We observe that in order to obtain the same result, adding $B\vec{e}_s$ needs to be accompanied by shifting $n$ to $n-1$ and $m$ to $m-1$. We have the same super $\hat{Z}$ after the action.

\subsection{The three variable knot invariant}

We use the results of the previous sections to define a simpler knot invariant. Specifically, the actions of the boundary torus of $Y_K$ on $\hat{Z}_{b,c}$ imply that infinitely many different $(b,c)$ pairs are related. Hence, we choose $\hat{Z}_{0,0}$ to be independent. Furthermore, they imply that  $\hat{Z}_{0,0}(M_K, y,z,n,m,q)$ is independent of values of $n,m \in \intg$. Using these properties of $\hat{Z}_{b,c}$, we define a three variable knot invariant.
$$
F_K (y,z,q) = F_K (y,z,q;\alpha_+)  + F_K (y,z,q;\alpha_-),
$$
\begin{align*}
F_K (y,z,q;\alpha_+) & : = \hat{Z}^{(\alpha_+)}_{0,0}(Y_K, y,z,n,m,q) \in  \intg + q^{\Delta}\intg[q^{-1},q]] [[y, z^{-1}]] \quad \text{in $\alpha_1$ chamber}   \\
F_K (y,z,q;\alpha_-) & : = \hat{Z}^{(\alpha_-)}_{0,0}(Y_K, y,z,n,m,q) \in  \intg + q^{\Delta}\intg[q^{-1},q]] [[z, y^{-1}]] \quad \text{in $\alpha_2 $ chamber}  \\ 
\end{align*}
where $\Delta \in \rational$ and $\intg[q^{-1},q]] [[y, z^{-1}]]$ denotes a vector space of Laurent power series in $y$ and $z^{-1}$ with coefficients in a Laurent power series ring $\intg[q^{-1},q]]$. Similarly for $\intg[q^{-1},q]] [[z, y^{-1}]]$. In each of the above good chambers, the results from Section 2 ensure that the power series is bounded below (i.e. well defined).
\newline

By the conjugation symmetry of the relative $Spin^c$ structures (15),
$$
F_K (1/y,1/z,q) = - F_K (y,z,q).
$$
Therefore, the general form of super $F_K$ is
\begin{equation}
F_K  (y,z,q) = c + \sum_{ \substack{m,n \in \intg_{\geq 0}^2 \\ (m,n) \neq (0,0)}} f_{m,n}(K; q) \lb \frac{y^m}{z^n} - \frac{z^n}{y^m} \rb \quad \in \intg + q^{\Delta}\intg[q^{-1},q]] [[y/z, (y/z)^{-1} ]].
\end{equation}
\begin{Rmk} The constant $c$ is finite in contrast to the diverging constant in (2) (cf. Remark 2.4).
\end{Rmk}
Compared to (38) in Appendix B, the summation in (38) is over odd integers and there is no chamber structure. Furthermore, we will see in Section 5.3 that the super $F_K$ of torus knots contains a classical series ($f_{m,n}(K; q) =$ constant for certain $(m,n)$ pairs) as well as $q$-dependent parts. To the best of the author's knowledge, there are no known examples of knots in which (38) contains such a classical series.

\section{Torus knots}

\subsection{Plumbing graphs}

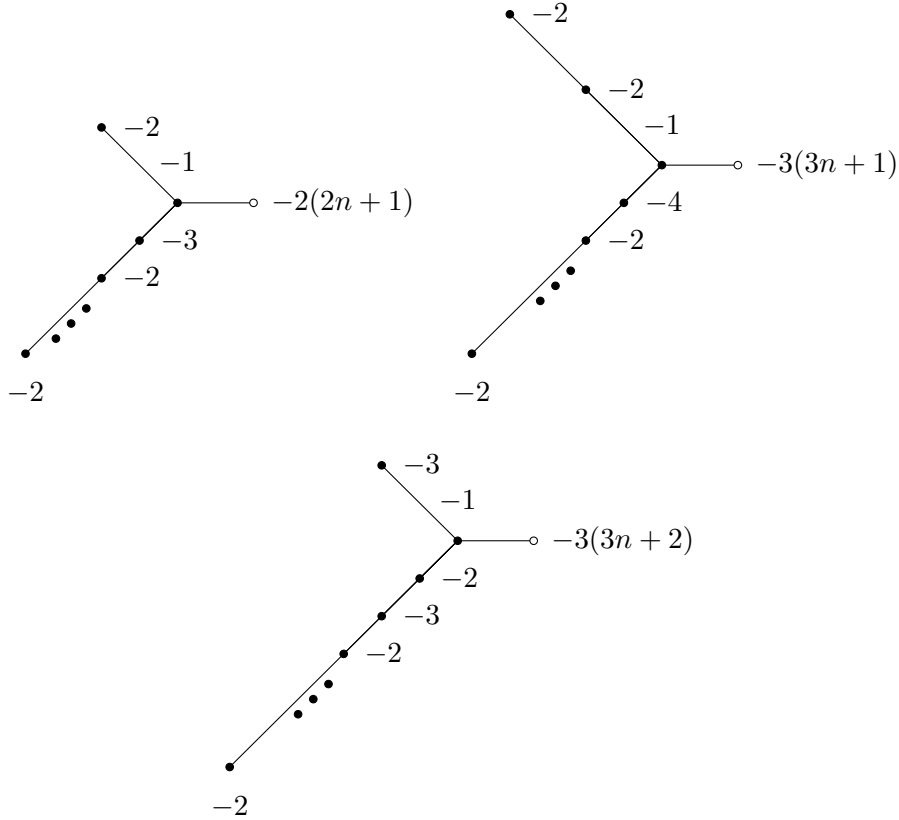
\begin{figure}[t!]
\begin{center}
\begin{tikzpicture}[scale=1]
\centering
\tikzstyle{every node}=[draw,shape=circle]

\draw (0,0)  node[circle,fill,inner sep=1pt,label=above:$-1$](){}; 
\draw (0,0) -- (1,0) node[circle,fill=white,inner sep=1pt,label=right:$-2(2n+1)$](){};
\draw (0,0)  -- (-1,1) node[circle,fill,inner sep=1pt,label=right:$-2$](-2){};

\draw (0,0) -- (-0.5,-0.5) node[circle,fill,inner sep=1pt,label=right:$-3$](-3){};
\draw (0,0)  -- (-1,-1) node[circle,fill,inner sep=1pt,label=right:$-2$](-2){};
\draw (0,0) -- (-2,-2) node[circle,fill,inner sep=1pt,label=below:$-2$](){};

\draw (-1.2,-1.4) node[circle,fill,inner sep=1pt,label=above:$$](){};
\draw (-1.4,-1.6) node[circle,fill,inner sep=1pt,label=above:$$](){};
\draw (-1.6,-1.8) node[circle,fill,inner sep=1pt,label=above:$$](){};

\end{tikzpicture}
\begin{tikzpicture}[scale=1]
\tikzstyle{every node}=[draw,shape=circle]

\draw (0,0) node[circle,fill,inner sep=1pt,label=above:$-1$](){};
\draw (0,0) -- (1,0) node[circle,fill=white,inner sep=1pt,label=right	:$-3(3n+1)$](){};
\draw (0,0)  -- (-1,1) node[circle,fill,inner sep=1pt,label=right:$-2$](-2){};
\draw (0,0)  -- (-2,2) node[circle,fill,inner sep=1pt,label=right:$-2$](-2){};

\draw (0,0) -- (-0.5,-0.5) node[circle,fill,inner sep=1pt,label=right:$-4$](-4){};
\draw (0,0)  -- (-1,-1) node[circle,fill,inner sep=1pt,label=right:$-2$](-2){};
\draw (0,0) -- (-2.5,-2.5) node[circle,fill,inner sep=1pt,label=below:$-2$](-2){};

\draw (-1.2,-1.4) node[circle,fill,inner sep=1pt,label=above:$$](){};
\draw (-1.4,-1.6) node[circle,fill,inner sep=1pt,label=above:$$](){};
\draw (-1.6,-1.8) node[circle,fill,inner sep=1pt,label=above:$$](){};

\end{tikzpicture}
\begin{tikzpicture}[scale=1]
\centering
\tikzstyle{every node}=[draw,shape=circle]

\draw (0,0) node[circle,fill,inner sep=1pt,label=above:$-1$](-1){}; 
\draw (0,0) -- (1,0) node[circle,fill=white,inner sep=1pt,label=right:$-3(3n+2)$](){};
\draw (0,0)  -- (-1,1) node[circle,fill,inner sep=1pt,label=right:$-3$](-3){};

\draw (0,0) -- (-0.5,-0.5) node[circle,fill,inner sep=1pt,label=right:$-2$](-2){};
\draw (0,0)  -- (-1,-1) node[circle,fill,inner sep=1pt,label=right:$-3$](-3){};
\draw (0,0)  -- (-1.5,-1.5) node[circle,fill,inner sep=1pt,label=right:$-2$](-2){};
\draw (0,0) -- (-3,-3) node[circle,fill,inner sep=1pt,label=below:$-2$](-2){};

\draw (-1.7,-1.9) node[circle,fill,inner sep=1pt,label=above:$$](){};
\draw (-1.9,-2.1) node[circle,fill,inner sep=1pt,label=above:$$](){};
\draw (-2.1,-2.3) node[circle,fill,inner sep=1pt,label=above:$$](){};

\end{tikzpicture}
\end{center}
\caption{Plumbing graphs of $T(2,2n+1)$ (left), $T(3,3n+1)$ (right) and, $T(3,3n+2)$ (bottom). The ellipsis indicates intermediate vertices with weight $-2$ along the legs. The total number of successive $-2$ vertices on the leg is $n-1$ for $T(2,2n+1), T(3,3n+1)$ and $T(3,3n+2)$.}
\end{figure}

We review the method for obtaining plumbing graphs of torus knots in \cite{GM} and then move on to finding good chambers for the knots. We next calculate examples of the super $F_K$.
\newline

We consider torus knots $T(s,t) \subset S^3$ where $\gcd(s,t) =1,\, 2 \leq s < t$. Torus knots are examples of algebraic knots. Their complements admit plumbing graph presentations. The graphs consist of one multivalent vertex having degree $3$ and weight $-1$ and three legs attached to the vertex. One of the legs has an open vertex of degree 1 called distinguished vertex representing a torus boundary of the knot complement. To find vertices and weights on the other legs, we solve 
$$
\frac{t^{\prime}}{t} + \frac{s^{\prime}}{s} = 1- \frac{1}{st}
$$
for unique integers $t^{\prime} \in (0,t)$ and $s^{\prime} \in (0,s)$ satisfying
$$
st^{\prime} \equiv -1 \, (\text{mod}\, t) \qquad ts^{\prime} \equiv -1 \, (\text{mod}\, s).
$$
Then we expand $-t/t^{\prime}$ and $-s/s^{\prime}$ as continued fractions, as described in Section 3.1. Each of them forms a leg with weights attached to the central vertex. The weight of the distinguished vertex is given by $-st$~\footnote{This value corresponds to $0$-framed torus knots.}. Examples of plumbing graphs are shown in Figure 4. 
\newline
\begin{Rmk} As in $sl(2)$ $F_K$ case~\cite{GM}, the super $F_K$ is applicable to torus knots in $\intg HS^3$. Figure 5 shows a method of obtaining a plumbing graph of the knots in $\intg HS^3$ from that of $S^3$. 
\end{Rmk}

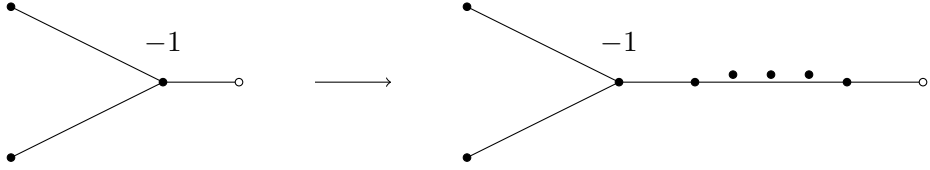
\begin{figure}
\begin{center}
\begin{tikzpicture}[scale=1]
\centering
\tikzstyle{every node}=[draw,shape=circle]

\draw (0,0) node[circle,fill,inner sep=1pt,label=above:$-1$](-1){} -- (1,0) node[circle,fill=white,inner sep=1pt,label](){};
\draw (0,0)  -- (-2,1) node[circle,fill,inner sep=1pt,label](){};
\draw (0,0) -- (-2,-1) node[circle,fill,inner sep=1pt,label](){};

\draw[->] (2,0) -- (3,0);

\draw (6,0)  -- (4,1) node[circle,fill,inner sep=1pt,label](){};
\draw (6,0) -- (4,-1) node[circle,fill,inner sep=1pt,label](){};

\draw  (6,0) node[circle,fill,inner sep=1pt,label=above:$-1$](-1){} -- (7,0) node[circle,fill,inner sep=1pt,label](){};
\draw  (7,0) --  (9,0) node[circle,fill,inner sep=1pt,label](){};
\draw  (9,0) --  (10,0) node[circle,fill=white,inner sep=1pt,label](){};

\draw (7.5,0.1) node[circle,fill,inner sep=1pt,label=above:](){};
\draw (8,0.1) node[circle,fill,inner sep=1pt,label=above:](){};
\draw (8.5,0.1) node[circle,fill,inner sep=1pt,label=above:](){};

\end{tikzpicture}

\end{center}
\caption{Changing the plumbing graph of $T(s,t) \subset S^3$ to that of $T(s,t) \subset \intg HS^3$. The graph without the distinguished vertex corresponds to a plumbing graph of $\intg HS^3$. The ellipsis indicates intermediate vertices.}
\end{figure}
	
\subsection{Chambers}

We find good chambers for infinite families of torus knots.

\begin{Prop} Let $v$ be the number of vertices of plumbing graphs of $T(2, 2n+1)$ and $T(3, 3n+ w), w=1,2$. Let $\alpha_{+}= (\alpha_1, \alpha_2 , \alpha_{v-1})$ and $\alpha_{-}$ be the good chambers for torus knots , where $\alpha_1$ corresponds to degree three vertex and the other two are associated with degree one vertices of their plumbing graphs. Their good chambers given by
$$
\alpha_{+} =  (1, 1, 1),\quad \alpha_{-} = -\alpha_{+},
$$
yield a well defined (Laurent) power series $f_{m,n}(q)$. 
\end{Prop}

\noindent \textit{Proof}. We describe a general strategy. The plumbing graphs of $T(s,t)$ consist of one degree 3 and three degree 1 vertices. The two regular vertices of the latter contribute to the integrand of (14). In order to find good chambers, we consider the regime of large powers of $q$ in the expansions (8) of the degree 3 vertex. Consequently, the prefactors of (8) have no effect. Then after denoting variables of the degree 3 vertex by $y_1$ and $z_1$, the relevant part of the integrand of $\hat{Z}$ for very large $r_1$ becomes
$$
\text{Integrand} \supset y_1^{\alpha_1 r_1} z_1^{-\alpha_1 r_1} 
\lb \sum\limits_{r_j } y_{j}^{\alpha_j r_j} + \sum\limits_{r_j } z_{j}^{- \alpha_j r_j} \rb
\lb \sum\limits_{r_{v-1} } y_{v-1}^{ \alpha_{v-1} r_{v-1}} + \sum\limits_{r_{v-1} } z_{v-1}^{- \alpha_{v-1} r_{v-1}} \rb
\prod\limits_{i=1}^{v} y_i^{l_{1,i}}z_i^{l_{2,i}}
$$
where $j$ and $v-1$ are the indices of the degree one vertices and $\alpha_1, \alpha_j, \alpha_{v-1} = \pm 1$. Using that $\vec{y}$ and $\vec{z}$ integrations in (14) extracts constant terms, $l_{1,i}$ and $l_{2,i}$ can be expressed in terms of $\vec{r}$ and $\vec{\alpha}$. This implies that we have a system of linear equations for $b=c=0$,
\begin{equation}
M \vec{n} = \vec{l_1}, \qquad  M \vec{m} = \vec{l_2} 
\end{equation}
where $n_v = m_v = 0$. There are several cases depending on the values of the right hand side of (28). In each case, we solve for $\vec{n}$ and $\vec{m}$ in terms of  $\vec{r}$ and $\vec{\alpha}$. Then substituting them into $q^{(\vec{n},B\vec{m})}$, we obtain $q^{f(B_{ij},\vec{r},\vec{\alpha})}$. From this we can determine whether or not good chambers exist.
\newline

It turns out that in most of cases of (28), the exponent $f(B_{ij},\vec{r},\vec{\alpha})$ is not bounded from below or not all $\alpha_i$'s appear in $f(B_{ij},\vec{r},\vec{\alpha})$. Hence, in the former case, the $q$ series is not convergent in $|q|<1$. In the latter case, complete chambers cannot be determined. There are two chambers that are good chambers. This corresponds to the case $\vec{l_1}= (-\alpha_1 r_1  ,-\alpha_2 r_2, -\alpha_3 r_3 , \ast )$ and $\vec{l_2}=  ( \alpha_1 r_1  ,0, 0 , \ast ) $.
\newline

\indent Applying the strategy to $T(2,2n+1)$ whose graph is depicted in Figure 4 we find that the exponent $f(B_{ij},\vec{r},\vec{\alpha})$ divides into three groups
\begin{multline*}
 \alpha_{1}^2 r_{1}^2 \lsb (2t)^2 |B_{11}| + \frac{(2t)^2}{4} |B_{22}| + \sum_{i=0}^{v-4} (t-2i-1)^2 | B_{i+3,i+3}| - (2t)^2 |B_{12}| - 2t(2t-2) |B_{13}| \right.\\ 
\left. -2 \sum_{i=0}^{v-4} (t-2i-1)(t-2i-3) | B_{i+3,i+4}| \rsb = 2t \alpha_{1}^2 r_{1}^2
\end{multline*} 
\begin{multline*}
\alpha_{1}\alpha_{2} r_{1}r_{2}  \lsb 2t^2 |B_{11}| + \frac{t(t+1)}{2} |B_{22}| + \frac{1}{2} \sum_{i=0}^{v-4} (t-2i-1)^2 | B_{i+3,i+3}| -t(2t+1)|B_{12}| \right.\\
\left.  -2t(t-1)|B_{13}| - \sum_{i=0}^{v-4} (t-2i-1)(t-2i-3) | B_{i+3,i+4}|  \rsb = t \alpha_{1}\alpha_{2}  r_{1}r_{2}
\end{multline*}
\begin{multline*}
 \alpha_{1}\alpha_{v-1} r_{1}r_{v-1}\lsb 4t|B_{11}|  + t|B_{22}| +  \sum_{i=0}^{v-4} (t-2i-1) | B_{i+3,i+3}| -4t|B_{12}|  \right.\\
\left.  -2(2t-1)|B_{13}| - 2 \sum_{i=0}^{v-4} (t-2i) | B_{i+2,i+3}|  \rsb = 2 \alpha_{1}\alpha_{v-1} r_{1}r_{v-1}
\end{multline*}
where $B_{11} = -1, B_{22} = -2,B_{33} = -3, B_{ii} = -2 (i>3),  B_{ij} = 1 (i\neq j)$ and $v=(t+5)/2$ are used. Combining them, we get
$$
f(B_{ij},\vec{r},\vec{\alpha}) = 2t \alpha_{1}^2 r_{1}^2  +  t \alpha_{1}\alpha_{2}  r_{1}r_{2} +  2 \alpha_{1}\alpha_{v-1} r_{1}r_{v-1},\qquad t=2n+1.
$$
This implies that 
\begin{equation}
\alpha_{1}= \alpha_{2} =  \alpha_{v-1} = 1\quad \text{or}\quad \alpha_{1}= \alpha_{2} =  \alpha_{v-1} = -1 
\end{equation}
ensures the boundedness of the super $F_K$ of $T(2,2n+1)$. Hence, they are good chambers.
\newline

In the case of $T(3,3n+1)$ in Figure 4, $f(B_{ij},\vec{r},\vec{\alpha})$ becomes
\begin{multline*}
 \alpha_{1}^2 r_{1}^2\lsb 9t^2 |B_{11}| + 4t^2 |B_{22}| + t^2 |B_{33}| + \sum_{i=0}^{v-5} (t-3i-1)^2 |B_{i+4,i+4}| - 12t^2 |B_{12}| - 6t(t-1)|B_{14}| \right. \\
\left. - 4t^2 |B_{23}| - 2 \sum_{i=1}^{v-5} (t-3i-1)(t-3i+2)  |B_{i+3,i+4}| \rsb = 3t \alpha_{1}^2 r_{1}^2 
\end{multline*}
\begin{multline*}
 \alpha_{1}\alpha_{3} r_{1}r_{3}\lsb 3t^2 |B_{11}| + \frac{2t(2t+1)}{3} |B_{22}| + \frac{t(t+2)}{3} |B_{33}| + \frac{1}{3} \sum_{i=0}^{v-5} (t-3i-1)^2 | B_{i+4,i+4}|\right.\\
\left. -t(4t+1)|B_{12}| -2t(t-1)|B_{14}| - \frac{1}{3} t(4t+5) |B_{23}| - \frac{2}{3} \sum_{i=1}^{v-4} (t-3i-1)(t-3i+2) | B_{i+3,i+4}|  \rsb 
\end{multline*}
$
= t \alpha_{1}\alpha_{3} r_{1}r_{3}
$
\begin{multline*}
 \alpha_{1}\alpha_{v-1} r_{1}r_{v-1}\lsb 9t|B_{11}|  + 4t|B_{22}| +  t|B_{33}| + \sum_{i=0}^{v-5} (t-3i-1) | B_{i+4,i+4}| \right.\\
\left. -12t |B_{12}|  -3(2t-1)|B_{14}| - 4t |B_{23}|  - \sum_{i=1}^{v-5} (2t-6i+1) | B_{i+3,i+4}|  \rsb =  3  \alpha_{1}\alpha_{v-1} r_{1}r_{v-1}
\end{multline*}
where $B_{11} = -1, B_{22} = -2,B_{33} = -2, B_{44} = -4, B_{ii} = -2 (i>4),  B_{ij} = 1 (i\neq j)$ and $v=(t+11)/3$ are used. Combining them, we get
$$
f(B_{ij},\vec{r},\vec{\alpha}) = 3t \alpha_{1}^2 r_{1}^2  +  t  \alpha_{1}\alpha_{2}  r_{1}r_{2} +  3 \alpha_{1}\alpha_{v-1} r_{1}r_{v-1},\qquad t=3n+1.
$$
We arrive at (29).
\newline

\indent In the case of $T(3,3n+2)$ in Figure 4, $f(B_{ij},\vec{r},\vec{\alpha})$ becomes
\begin{multline*}
 \alpha_{1}^2 r_{1}^2\lsb 9t^2 |B_{11}| + t^2 |B_{22}| + (2t-1)^2 |B_{33}| + \sum_{i=0}^{v-5} (t-3i-1)^2 |B_{i+4,i+4}| - 6t^2 |B_{12}| - 6t(2t-1)|B_{13}| \right.\\
\left. -2(t-2)(2t-1) |B_{34}| - 2 \sum_{i=1}^{v-5} (t-3i-2)(t-3i+1)  |B_{i+3,i+4}| \rsb = 3t \alpha_{1}^2 r_{1}^2 
\end{multline*}
\begin{multline*}
 \alpha_{1}\alpha_{2} r_{1}r_{2}\lsb 3t^2 |B_{11}| + \frac{t(t+1)}{3} |B_{22}| + \frac{(2t-1)^2}{3} |B_{33}| + \frac{1}{3} \sum_{i=0}^{v-5} (t-3i-2)^2 | B_{i+4,i+4}|\right.\\
\left. -t(2t+1)|B_{12}| -2t(2t-1)|B_{13}| - \frac{2}{3}(2t-1)(t-2) |B_{34}| - \frac{2}{3} \sum_{i=1}^{v-4} (t-3i-2)(t-3i+1) | B_{i+3,i+4}| \rsb \\
\end{multline*}
$
= t \alpha_{1}\alpha_{2} r_{1}r_{2}
$
\begin{multline*}
 \alpha_{1}\alpha_{v-1} r_{1}r_{v-1}\lsb 9t|B_{11}|  + t|B_{22}| +  2(2t-1) |B_{33}| + \sum_{i=0}^{v-5} (t-3i-2) | B_{i+4,i+4}| \right.\\
\left. -6t |B_{12}|  -3(4t-1)|B_{13}| - (4t-5) |B_{34}|  - \sum_{i=1}^{v-5} (2t-6i-1) | B_{i+3,i+4}|  \rsb =  3  \alpha_{1}\alpha_{v-1} r_{1}r_{v-1}
\end{multline*}
where $B_{11} = -1, B_{22} = -3, B_{33} = -2, B_{44} = -3, B_{ii} = -2 (i>4),  B_{ij} = 1 (i\neq j)$ and $v=(t+10)/3$ are used. We obtain the same result as (30).
\newline
\begin{Rmk}  In Proposition 5.2, it is expected that the values of $\alpha_{+,2}$ and $\alpha_{+,v-1}$ are the same because the corresponding vertices have degree one and the degree two vertices between the degree one and three vertices have no effect.
\end{Rmk}

\begin{Conj}
Proposition 5.2 holds for all torus knots $T(s,t) \subset S^3\, (gcd(s,t)=1)$.
\end{Conj}

\subsection{Examples}

We apply the results of the previous sections to calculate examples of super $F_K$. Additional examples are recorded in Appendix A.
\newline

\noindent $\underline{K=T(2,3)}$ Using Figure 4, (14) and Proposition 5.2, we obtain
$$
F_{T(2,3)}(y,z,q)  =  1 + \sum_{i=2}^{\infty} \lb y^i +  \frac{1}{z^i} \rb - \sum_{i=2}^{\infty} \lb \frac{1}{y^i} + z^i \rb +q \left(\frac{y^2}{z^3}+\frac{y^3}{z^2}  - \frac{z^2}{y^3}-\frac{z^3}{y^2} \right)
$$
$$
+q^2 \left(\frac{y^3}{z^4}+\frac{y^4}{z^3} -\frac{z^3}{y^4}-\frac{z^4}{y^3}\right) +q^5 \left(-\frac{y^5}{z^6}-\frac{y^6}{z^5} +\frac{z^5}{y^6}+\frac{z^6}{y^5}\right)+q^7 \left(-\frac{y^6}{z^7}-\frac{y^7}{z^6} +\frac{z^6}{y^7}+\frac{z^7}{y^6}\right) 
$$
\begin{equation}
 +q^{12}\left(\frac{y^8}{z^9}+\frac{y^9}{z^8} - \frac{z^8}{y^9}-\frac{z^9}{y^8}\right)  +q^{15} \left(\frac{y^9}{z^{10}}+\frac{y^{10}}{z^9} -\frac{z^9}{y^{10}}-\frac{z^{10}}{y^9}\right) + \cdots
\end{equation}
%
\noindent We find that (31) splits into $q$ independent and dependent parts. The former is a new feature in the super $F_K$, which is absent in the $F_K$ associated with $sl(2)$ (cf. (39) in Appendix B). The role of the former will be described in Section 6. It can be expressed in terms of the unknot (9).
$$
\sum_{i=2}^{\infty} \lb y^i +  \frac{1}{z^i} \rb - \sum_{i=2}^{\infty} \lb \frac{1}{y^i} + z^i \rb = \frac{y-z}{(1-y)(1-z)}\bigg|_{\alpha_+}+  \frac{y-z}{(1-y)(1-z)}\bigg|_{\alpha_-}   - \lb y - \frac{1}{y} \rb + \lb z - \frac{1}{z} \rb  -2.
$$
The latter can be cast into (26) given by
$$
f_{m,n}(K; q) = \epsilon_{m,n} q^{\frac{mn}{6}}, \qquad n=m+1,
$$
$$
\epsilon_{m,n}(K)  = \begin{cases}
+1,\, r^{+-}_m \equiv 3 \quad\& \quad  r^{+-}_n \equiv 4 \quad\& \quad r^{-+}_m \equiv 1 \quad\& \quad r^{-+}_n \equiv 2 \quad \text{mod} \, 6\\
+1,\, r^{+-}_m \equiv 4 \quad\& \quad  r^{+-}_n \equiv 5 \quad\& \quad   r^{-+}_m \equiv 2\quad\& \quad r^{-+}_n \equiv 3  \quad \text{mod} \, 6 \\
-1,\, r^{++}_m \equiv 0 \quad\& \quad  r^{++}_n \equiv 1 \quad\& \quad   r^{--}_m \equiv 4 \quad\& \quad  r^{--}_n \equiv 5  \quad \text{mod} \, 6 \\
-1,\, r^{++}_m \equiv 1 \quad\& \quad  r^{++}_n \equiv 2 \quad\& \quad   r^{--}_m \equiv 5 \quad\& \quad  r^{--}_n \equiv 0   \quad \text{mod} \, 6 \\
0, \quad \text{otherwise}
\end{cases}
$$
where $r^{++}_m = m - (6+2+3), r^{--}_m = m - (6-2-3), r^{+-}_m = m - (6+2-3), r^{-+}_m = m - (6-2+3)$ and $r_n$'s are obtained by replacing $m$ by $n$ from $r_m$'s. 
\newline

\noindent $\underline{K=T(2,5)}$ We obtain
$$
F_{T(2,5)}(y,z,q) = 1 + \sum_{\substack{i=2 \\ i\neq 3}}^{\infty} \lb y^i + \frac{1}{z^i} \rb  - \sum_{\substack{i=2 \\ i\neq 3}}^{\infty} \lb \frac{1}{y^i} + z^i \rb+q \left(\frac{y^2}{z^5}+\frac{y^5}{z^2}-\frac{z^2}{y^5}-\frac{z^5}{y^2}\right)
$$
$$
+ q^2 \left(\frac{y^4}{z^5}+\frac{y^5}{z^4}-\frac{z^4}{y^5}-\frac{z^5}{y^4}\right)+q^3\left(\frac{y^5}{z^6}+\frac{y^6}{z^5}-\frac{z^5}{y^6}-\frac{z^6}{y^5}\right)+q^4 \left(\frac{y^5}{z^8}+\frac{y^8}{z^5}-\frac{z^5}{y^8}-\frac{z^8}{y^5}\right)
$$
\begin{equation}
+q^7 \left(-\frac{y^7}{z^{10}}-\frac{y^{10}}{z^7} +\frac{z^7}{y^{10}}+\frac{z^{10}}{y^7}\right) +q^9 \left(-\frac{y^9}{z^{10}}-\frac{y^{10}}{z^9} +\frac{z^9}{y^{10}}+\frac{z^{10}}{y^9}\right) + \cdots
\end{equation}
\noindent We again find that (32) splits into $q$ independent and dependent parts. The former can also be expressed in terms of (9).
\begin{align*}
\sum_{\substack{i=2 \\ i\neq 3}}^{\infty} \lb y^i + \frac{1}{z^i} \rb  - \sum_{\substack{i=2 \\ i\neq 3}}^{\infty} \lb \frac{1}{y^i} + z^i \rb & = \frac{y-z}{(1-y)(1-z)}\bigg|_{\alpha_+}+  \frac{y-z}{(1-y)(1-z)}\bigg|_{\alpha_-}   - \lb y^3 + y -  \frac{1}{y} - \frac{1}{y^3} \rb\\
&  +  z + z^3  - \frac{1}{z} - \frac{1}{z^3} -2.
\end{align*}
The latter can be cast into (27) given by
\begin{equation}
\epsilon_{m,n} q^{\frac{m(m+g(m,n))}{10}}\lb \frac{y^m}{z^{m+g(m,n)}} + \frac{y^{m+g(m,n)}}{z^{m}} -  \frac{z^m}{y^{m+g(m,n)}} - \frac{z^{m+g(m,n)}}{y^{m}}\rb,
\end{equation}
$$
\epsilon_{m,n}(K) = \begin{cases}
+1,\, r^{+-}_m \equiv 5 \quad\& \quad  r^{+-}_n \equiv 8 \quad\& \quad   r^{-+}_m \equiv 9 \quad\& \quad r^{-+}_n \equiv 2 \quad \text{mod} \, 10\\
+1,\, r^{+-}_m \equiv 7 \quad\& \quad  r^{+-}_n \equiv 8 \quad\& \quad   r^{-+}_m \equiv 1 \quad\& \quad r^{-+}_n \equiv 2  \quad \text{mod} \, 10 \\
+1,\, r^{+-}_m \equiv 8 \quad\& \quad  r^{+-}_n \equiv 9 \quad\& \quad   r^{-+}_m \equiv 2 \quad\& \quad r^{-+}_n \equiv 3  \quad \text{mod} \, 10 \\
+1,\, r^{+-}_m \equiv 8 \quad\& \quad  r^{+-}_n \equiv 1 \quad\& \quad   r^{-+}_m \equiv 2 \quad\& \quad r^{-+}_n \equiv 5  \quad \text{mod} \, 10 \\
-1,\, r^{++}_m \equiv 0 \quad\& \quad  r^{++}_n \equiv 3 \quad\& \quad   r^{--}_m \equiv 4 \quad\& \quad  r^{--}_n \equiv 7  \quad \text{mod} \, 10 \\
-1,\, r^{++}_m \equiv 2 \quad\& \quad  r^{++}_n \equiv 3 \quad\& \quad   r^{--}_m \equiv 6 \quad\& \quad  r^{--}_n \equiv 7   \quad \text{mod} \,10 \\
-1,\, r^{++}_m \equiv 3 \quad\& \quad  r^{++}_n \equiv 4 \quad\& \quad   r^{--}_m \equiv 7 \quad\& \quad  r^{--}_n \equiv 8   \quad \text{mod} \,10 \\
-1,\, r^{++}_m \equiv 3 \quad\& \quad  r^{++}_n \equiv 6 \quad\& \quad   r^{--}_m \equiv 7 \quad\& \quad  r^{--}_n \equiv 0   \quad \text{mod} \,10 \\
0, \quad \text{otherwise}
\end{cases}
$$
where $r^{++}_m = m - (10+2+5), r^{--}_m = m - (10-2-5), r^{+-}_m = m - (10+2-5), r^{-+}_m = m - (10-2+5)$ and $r_n$'s can be obtained by replacing $m$ by $n$.
$$
g(m,n) = \begin{cases}
3,\, r^{+-}_m \equiv 5 \quad\& \quad  r^{+-}_n \equiv 8 \quad\& \quad   r^{-+}_m \equiv 9 \quad\& \quad r^{-+}_n \equiv 2 \quad \text{mod} \, 10\\
3,\, r^{+-}_m \equiv 8 \quad\& \quad  r^{+-}_n \equiv 1 \quad\& \quad   r^{-+}_m \equiv 2 \quad\& \quad r^{-+}_n \equiv 5  \quad \text{mod} \, 10 \\
3,\, r^{++}_m \equiv 0 \quad\& \quad  r^{++}_n \equiv 3 \quad\& \quad   r^{--}_m \equiv 4 \quad\& \quad  r^{--}_n \equiv 7  \quad \text{mod} \, 10 \\
3,\, r^{++}_m \equiv 3 \quad\& \quad  r^{++}_n \equiv 6 \quad\& \quad   r^{--}_m \equiv 7 \quad\& \quad  r^{--}_n \equiv 0   \quad \text{mod} \,10 \\
1,\, r^{+-}_m \equiv 7 \quad\& \quad  r^{+-}_n \equiv 8 \quad\& \quad   r^{-+}_m \equiv 1 \quad\& \quad r^{-+}_n \equiv 2  \quad \text{mod} \, 10 \\
1,\,  r^{+-}_m \equiv 8 \quad\& \quad  r^{+-}_n \equiv 9 \quad\& \quad   r^{-+}_m \equiv 2 \quad\& \quad r^{-+}_n \equiv 3  \quad \text{mod} \, 10 \\
1,\, r^{++}_m \equiv 2 \quad\& \quad  r^{++}_n \equiv 3 \quad\& \quad   r^{--}_m \equiv 6 \quad\& \quad  r^{--}_n \equiv 7   \quad \text{mod} \,10 \\
1,\, r^{++}_m \equiv 3 \quad\& \quad  r^{++}_n \equiv 4 \quad\& \quad   r^{--}_m \equiv 7 \quad\& \quad  r^{--}_n \equiv 8   \quad \text{mod} \,10 \\
0, \quad \text{otherwise}
\end{cases}
$$
\newline

\noindent $\underline{K=T(2,7)}$ We obtain
$$
F_{T(2,7)}(y,z,q) = 1 + \sum_{\substack{i=2 \\ i\neq 3,5}}^{\infty} \lb y^i + \frac{1}{z^i} \rb  - \sum_{\substack{i=2 \\ i\neq 3,5}}^{\infty} \lb \frac{1}{y^i} + z^i \rb + q \left(\frac{y^2}{z^7}+\frac{y^7}{z^2}-\frac{z^2}{y^7}-\frac{z^7}{y^2}\right)
$$
$$
+q^2\left(\frac{y^4}{z^7}+\frac{y^7}{z^4}-\frac{z^4}{y^7}-\frac{z^7}{y^4}\right)+q^3\left(\frac{y^6}{z^7}+\frac{y^7}{z^6}-\frac{z^6}{y^7}-\frac{z^7}{y^6}\right)+q^4 \left(\frac{y^7}{z^8}+\frac{y^8}{z^7}-\frac{z^7}{y^8}-\frac{z^8}{y^7}\right)
$$	
\begin{equation}
+q^5\left(\frac{y^7}{z^{10}}+\frac{y^{10}}{z^7}-\frac{z^7}{y^{10}}-\frac{z^{10}}{y^7}\right)+q^6\left(\frac{y^7}{z^{12}}+\frac{y^{12}}{z^7}-\frac{z^7}{y^{12}}-\frac{z^{12}}{y^7}\right) + q^9 \left(-\frac{y^9}{z^{14}}-\frac{y^{14}}{z^9}+\frac{z^9}{y^{14}}+\frac{z^{14}}{y^9}\right) +\cdots
\end{equation}
\noindent The $q$ independent terms can be expressed in terms of (9).
\begin{align*}
\sum_{\substack{i=2 \\ i\neq 3,5}}^{\infty} \lb y^i + \frac{1}{z^i} \rb  - \sum_{\substack{i=2 \\ i\neq 3,5}}^{\infty} \lb \frac{1}{y^i} + z^i \rb & = \frac{y-z}{(1-y)(1-z)}\bigg|_{\alpha_+}+  \frac{y-z}{(1-y)(1-z)}\bigg|_{\alpha_-}   - \lb y + y^3 + y^5 -  \frac{1}{y} \right. \\
 & \left.  -\frac{1}{y^3} -\frac{1}{y^5} \rb  + z + z^3  + z^5 -\frac{1}{z}  - \frac{1}{z^3} - \frac{1}{z^5}   -2
\end{align*}
The $q$ dependent terms of (34) can be cast into (33), where $g(m,n)$ and $\epsilon_{m,n}$ are in Appendix A.\\
\newline We observe that the super $F_K$ of the above knots is symmetric under the exchange of $m$ and $n$; thus $f_{n,m}= f_{m,n}$. Furthermore, the splitting structure of the super $F_K$ of the torus knots is similar to that of the multivariable knot polynomial associated with $sl(2|1)$ defined in \cite{GP1}. A difference is that the knot dependent part in the latter is a polynomial.
\newline

\noindent \underline{An Algorithm} We present a simple algorithm for finding the sign function $\epsilon_{m,n}$ and the exponent shift function $g(m,n)$ of the super $F_K$ for $T(2,2l+1), l\geq 2$~\footnote{The trefoil ($l=1$) is a special case of the algorithm.}. We recall that the $q$ dependent part is
$$
\epsilon_{m,n}  q^{\frac{m(m+g(m,n))}{2(2l+1)}}\lb \frac{y^m}{z^{m+g(m,n)}} + \frac{y^{m+g(m,n)}}{z^{m}} -  \frac{z^m}{y^{m+g(m,n)}} - \frac{z^{m+g(m,n)}}{y^{m}}\rb.
$$

\begin{enumerate}
	\item The shift function $g(m,n)$ takes values in $\lac 1,3,5,\cdots, 2l-1 \rac$. We let 
	  	\begin{align*}
			  r_{m}^{++} & = m - (2(2l+1)+2+2l+1)\\
				r_{m}^{--} & = m - (2(2l+1)-2-2l-1)\\
				r_{m}^{+-} & = m - (2(2l+1)+2-2l-1)\\
				r_{m}^{-+} & = m - (2(2l+1)-2+2l+1)\\			
			\end{align*}
	  The terms $r_{n}^{\pm\pm}$ are defined similarly.
				
	\item For the $-1$ case of $\epsilon_{m,n}$: Start with $max(g(m,n))$ and set $r_{m}^{++} \equiv 0$ and $r_{m}^{--} \equiv 4$.
	
	\item Set $r_{n}^{++} \equiv max(g(m,n))$ and $r_{n}^{--} \equiv 4 + max(g(m,n))$. We denote this pair by $p=(max(g(m,n)), 4+ max(g(m,n)))$.
	
	\item Move on to $max(g(m,n))-2$ and set $(r_{n}^{++}, r_{n}^{--}) \equiv p$. Next set $(r_{m}^{++} , r_{m}^{--}) \equiv p - (max(g(m,n))-2)$.
	
	
	\item Iterate Step 4 until $g(m,n)=1$ is completed ($r_{n}^{++}, r_{n}^{--}$ remain equal to $p$ during the iteration).

	\item Start again with $g(m,n)=1$. Set $( r_{m}^{++} , r_{m}^{--})\equiv p$ and $( r_{n}^{++} , r_{n}^{--}) \equiv p+ 1$.
	
	\item Move on to $g(m,n)=3$ and set $(r_{m}^{++}, r_{m}^{--}) \equiv p$. Next set $(r_{n}^{++} , r_{n}^{--}) \equiv p+ 3$.
	
	\item Iterate step 7 until $max(g(m,n))$ is reached.
	
	\item For the $+1$ case of $\epsilon_{m,n}$: Start with $max(g(m,n))$ and set $r_{m}^{+-} \equiv 2l+1$ and $r_{m}^{-+} \equiv 2l+5$.
	
	\item Set $r_{n}^{+-} \equiv 2l+1 + max(g(m,n))$ and $r_{n}^{-+} \equiv 2$.  We denote this pair by $s=(2l+1 + max(g(m,n)), 2)$.
	
	\item Move on to $max(g(m,n))-2$ and set $(r_{n}^{+-} , r_{n}^{-+}) \equiv s$. Next set $r_{m}^{+-} \equiv r_{n}^{+-} - ( max(g(m,n))-2)$ and $r_{m}^{-+} \equiv r_{m}^{+-} + 4$.
	
	\item Iterate Step 11 until $g(m,n)=1$ is completed ($r_{n}^{+-}, r_{n}^{-+}$ remain equal to $s$ during the iteration).
	
	\item Start again with $g(m,n)=1$. Set $( r_{m}^{+-} , r_{m}^{-+}) \equiv s$ and $(r_{n}^{+-}, r_{n}^{-+}) \equiv s + 1$.
	
	\item Move on to $g(m,n)=3$ and set $( r_{m}^{+-}, r_{m}^{-+} ) \equiv s$. Next set $r_{n}^{+-}\equiv 1$ and $r_{n}^{-+}\equiv 5$.
	
	\item Move on to $g(m,n)=5$ and set $( r_{m}^{+-} , r_{m}^{-+}) \equiv s$. Add 2 to the previous $r_{n}^{+-}$ and $r_{n}^{-+}$  to obtain $r_{n}^{+-}$ and $r_{n}^{-+}$ associated with $g=5$.
	
	\item Iterate Step 15 until $max(g(m,n))$ is reached ($r_{m}^{+-}, r_{m}^{-+}$ remain equal to $s$ during the iteration). The last $r_{n}^{+-}$ and $r_{n}^{-+}$ are $\equiv 2l-3$ and $2l+1$, respectively.
	
	\item Collect the cases for each value of $g(m,n)$. 

\end{enumerate}

\begin{Rmk} The case $\epsilon_{m,n}=-1$ consists of $(r_{m}^{++}$,  $r_{m}^{--}$ , $r_{n}^{++}$,  $r_{n}^{--})$ for each condition whereas $\epsilon_{m,n}=+1$ consists of $(r_{m}^{+-}$,  $r_{m}^{-+}$ , $r_{n}^{+-}$,  $r_{n}^{-+})$.
\end{Rmk}
\begin{Rmk} All the above conditions are over  $\intg/(2(2l+1))\intg$.
\end{Rmk}
\begin{Rmk} The shift function $g(m,n)$ is automatically generated by the algorithm.
\end{Rmk}
\begin{Rmk} We note that $2|g(m,n)|$ cases exist for each sign $\epsilon_{m,n}= \pm 1$.
\end{Rmk}

\noindent We apply the algorithm to $T(2,5)$ and $T(2,7)$. For $T(2,5)$,

\begin{enumerate}
	\item The exponent shift function is $g=\lac 1,3 \rac$ 
	 	\begin{align*}
			  r_{m}^{++} & = m - (10+2+5)\\
				r_{m}^{--} & = m - (10-2-5)\\
				r_{m}^{+-} & = m - (10+2-5)\\
				r_{m}^{-+} & = m - (10-2+5)\\			
			\end{align*}
			
		\item For $g=3$, $r_{m}^{++} \equiv 0$ and $r_{m}^{--}\equiv 4$.
		
		\item $r_{n}^{++} \equiv 3,  r_{m}^{--}\equiv 4+3 =7$. Then $p=(3,7)$. These, together with the above, are conditions for $g=3$
		
		\item For $g=1$, $r_{n}^{++} \equiv 3,  r_{n}^{--}\equiv 4+3 =7$, $r_{m}^{++} \equiv 3-1=2$ and $r_{m}^{--}\equiv 7-1=6$.
		
		\item This is the same as Step 4
		
		\item For another condition where $g=1$, $r_{m}^{++} \equiv 3, r_{m}^{--}\equiv 7$ and $r_{n}^{++} \equiv 4,  r_{n}^{--}\equiv 8$.
		
		\item For another condition where $g=3$, $r_{m}^{++} \equiv 3, r_{m}^{--}\equiv 7$ and $r_{n}^{++} \equiv 6,  r_{n}^{--}\equiv 0$.
		
		\item This is the same as Step 7
		
		\item For $g=3$, $r_{m}^{+-} \equiv 5, r_{m}^{-+}\equiv 9$.
		
		\item $r_{n}^{+-} \equiv 5+3=8, r_{n}^{-+}\equiv 2$. Then $s=(8,2)$. These, together with the above, are conditions for $g=3$
		
		\item For $g=1$, $r_{n}^{+-} \equiv 8, r_{n}^{-+}\equiv 2$ and $r_{m}^{+-} \equiv 8-1=7, r_{m}^{-+}\equiv 7+4=11\equiv 1$.
		
		\item Same as Step 11
		
		\item For another condition where $g=1$, $r_{m}^{+-} \equiv 8, r_{m}^{-+}\equiv 2$ and $r_{n}^{+-} \equiv 9, r_{n}^{-+}\equiv 3$.
		
		\item For another condition where $g=3$, $r_{m}^{+-} \equiv 8, r_{m}^{-+}\equiv 2$ and $r_{n}^{+-} \equiv 1, r_{n}^{-+}\equiv 5$.
	
\end{enumerate}
From Step 17, we arrive at $\epsilon_{m,n}$ and $g(m,n)$ in (33); Steps 15 and 16 are irrelevant. The application to $T(2,7)$ is recorded in Appendix A.
\newline

\noindent \underline{Reduction to $sl(2)$ } We provide an algorithm to extract $F^{sl(2,\complex)}_K (x,q)$ from $F_K (y,z,q)$ for $T(2,2l+1)$. 
\newline

\noindent An Algorithm:
\begin{enumerate}

  \item Let $\tilde{F}_K = q^l F^{+}_K (y,z,q)$, where $F^{+}_K (y,z,q)$ denotes the $q$ dependent $(y^b/z^c + y^c/z^b)$ part of $F_K (y,z,q)$.
	
	\item Pick the leading order (LO) term (i.e. the term with lowest $q$ power) in $\tilde{F}_K$. Skip the next $2(l-1)$ terms and then add the next two consecutive terms (i.e. the $2l$-th and $(2l+1)$-th terms) to the LO term.
	
 \item Skip the next $2(l-1)$ terms and then add the next two consecutive terms to the result of Step 2.
	
	\item Iterate Step 3) until a sufficient number of terms is reached.

	\item Eliminate terms of the form $y^r /z^w$, where $r>w$.
	
	\item Replace $y^r /z^w$ by $x^{\frac{r+w}{2}}$ and denote the resulting expression by $\tilde{F}_{K}^{+} (x,q)$.
	
	\item $F_{K}^{+}(x,q) = - q^l x^{\frac{2l-1}{2}} + \tilde{F}_{K}^{+} (x,q)$
	
	\item Perform Weyl symmetrization to obtain $F_{K}(x,q) = F_{K}^{+}(x,q) - F_{K}^{+}(1/x,q)$

\end{enumerate}
\begin{Rmk} In case of $T(2,3), (l=1)$, every term in $F^{+}_K (y,z,q)$ contribute to its $F_{K}(x,q)$.
\end{Rmk}

\noindent We apply the algorithm to $T(2,5)$.

\begin{enumerate}

\item From (32) we get \begin{align*}
\tilde{F}_K & = q^3 \left(\frac{y^2}{z^5}+\frac{y^5}{z^2}\right)+q^4 \left(\frac{y^4}{z^5}+\frac{y^5}{z^4}\right)+q^5\left(\frac{y^5}{z^6}+\frac{y^6}{z^5}\right)+q^6 \left(\frac{y^5}{z^8}+\frac{y^8}{z^5}\right)+ q^9 \left(-\frac{y^7}{z^{10}}-\frac{y^{10}}{z^7}\right)\\
& +q^{11} \left(-\frac{y^9}{z^{10}}-\frac{y^{10}}{z^9}\right) +q^{13}\left(-\frac{y^{10}}{z^{11}}-\frac{y^{11}}{z^{10}}\right)+q^{15}\left(-\frac{y^{10}}{z^{13}}-\frac{y^{13}}{z^{10}}\right)+q^{20}\left(\frac{y^{12}}{z^{15}}+\frac{y^{15}}{z^{12}}\right)\\
& +q^{23}\left(\frac{y^{14}}{z^{15}}+\frac{y^{15}}{z^{14}}\right) +q^{26} \left(\frac{y^{15}}{z^{16}}+\frac{y^{16}}{z^{15}}\right) +q^{29}\left(\frac{y^{15}}{z^{18}}+\frac{y^{18}}{z^{15}}\right)+q^{36}\left(-\frac{y^{17}}{z^{20}}-\frac{y^{20}}{z^{17}}\right) + \cdots
\end{align*}

\item The LO term is $q^3 \left(\frac{y^2}{z^5}+\frac{y^5}{z^2}\right)$. Add $q^6 \left(\frac{y^5}{z^8}+\frac{y^8}{z^5}\right)+q^9 \left(-\frac{y^7}{z^{10}}-\frac{y^{10}}{z^7}\right)$: 
$$
q^3 \left(\frac{y^2}{z^5}+\frac{y^5}{z^2}\right) + q^6 \left(\frac{y^5}{z^8}+\frac{y^8}{z^5}\right)+ q^9 \left(-\frac{y^7}{z^{10}}-\frac{y^{10}}{z^7}\right) 
$$

\item  $
q^3 \left(\frac{y^2}{z^5}+\frac{y^5}{z^2}\right) + q^6 \left(\frac{y^5}{z^8}+\frac{y^8}{z^5}\right)+ q^9 \left(-\frac{y^7}{z^{10}}-\frac{y^{10}}{z^7}\right) +q^{15}\left(-\frac{y^{10}}{z^{13}}-\frac{y^{13}}{z^{10}}\right)+q^{20}\left(\frac{y^{12}}{z^{15}}+\frac{y^{15}}{z^{12}}\right) 
$

\item  $
q^3 \left(\frac{y^2}{z^5}+\frac{y^5}{z^2}\right) + q^6 \left(\frac{y^5}{z^8}+\frac{y^8}{z^5}\right)+ q^9 \left(-\frac{y^7}{z^{10}}-\frac{y^{10}}{z^7}\right) +q^{15}\left(-\frac{y^{10}}{z^{13}}-\frac{y^{13}}{z^{10}}\right)+q^{20}\left(\frac{y^{12}}{z^{15}}+\frac{y^{15}}{z^{12}}\right) +q^{29}\left(\frac{y^{15}}{z^{18}}+\frac{y^{18}}{z^{15}}\right)+q^{36}\left(-\frac{y^{17}}{z^{20}}-\frac{y^{20}}{z^{17}}\right)   + \cdots
$

\item $
q^3 \frac{y^2}{z^5} + q^6 \frac{y^5}{z^8} - q^9 \frac{y^7}{z^{10}} - q^{15} \frac{y^{10}}{z^{13}} +q^{20} \frac{y^{12}}{z^{15}} +q^{29} \frac{y^{15}}{z^{18}} - q^{36}\frac{y^{17}}{z^{20}} + \cdots
$

\item $
\hat{F}_{K}^{+} (x,q) = q^3 x^{\frac{7}{2}} + q^6 x^{\frac{13}{2}} - q^9 x^{\frac{17}{2}} - q^{15}x^{\frac{23}{2}} + q^{20}x^{\frac{27}{2}} +q^{29} x^{\frac{33}{2}} - q^{36} x^{\frac{37}{2}} + \cdots
$

\item $F_{K}^{+}(x,q) = - q^2 x^{\frac{3}{2}} + q^3 x^{\frac{7}{2}} + q^6 x^{\frac{13}{2}} - q^9 x^{\frac{17}{2}} - q^{15}x^{\frac{23}{2}} + q^{20}x^{\frac{27}{2}} +  q^{29} x^{\frac{33}{2}} - q^{36} x^{\frac{37}{2}}  + \cdots$

\item $F_{K} (x,q) = - q^2 \lb x^{\frac{3}{2}} - x^{\frac{-3}{2}}\rb + q^3 \lb x^{\frac{7}{2}} - x^{\frac{-7}{2}} \rb + q^6 \lb x^{\frac{13}{2}} - x^{\frac{-13}{2}} \rb - q^9 \lb x^{\frac{17}{2}} -  x^{\frac{-17}{2}}\rb - q^{15} \lb x^{\frac{23}{2}} - x^{\frac{-23}{2}} \rb + q^{20} \lb x^{\frac{27}{2}} - x^{\frac{-27}{2}} \rb  +  q^{29} \lb x^{\frac{33}{2}} -  x^{\frac{-33}{2}} \rb -  q^{36}\lb x^{\frac{37}{2}} - x^{\frac{-37}{2}} \rb + \cdots$
	
\end{enumerate}
This agrees with (39) after reinstating the overall factor $q^{(s - 1)(t - 1)/2  - (st - s - t)^2 /4st}$.

\subsection{Mirror knots}

Polynomial invariants of a knot $K$--for example, the colored Jones polynomials and the HOMFLY-PT polynomials--behave simply under mirror reflection of knots $K^{\ast}$. To obtain the invariants of $K^{\ast}$, we send $q \mapsto q^{-1}$. In the case of the $F_K$ associated with the Lie algebra (35) in Appendix B, sending $q \mapsto q^{-1}$ is valid if the coefficient functions of $q$ are (Laurent) polynomials~\cite{GM}. However, there are knots such as $5_2$ whose coefficient functions are Laurent power series. For these knots, the above simple map is invalid. Examples of such knots were analyzed in \cite{P3}. This also applies to the super $F_K$. It behaves simply under mirror reflection if the coefficient functions $f_{m,n}(q)$ are (Laurent) polynomials. Examples of such knots are algebraic knots in $S^3$, which in turn contain torus knots. Specifically, under mirror reflection, the chambers in (27) are exchanged,
\newline
%
$$
 F_K (y,z,q;\alpha_1) \mapsto  F_{K^{\ast}}  (y^{-1} ,z^{-1} , q^{-1}; \alpha_1) =   F_{K^{\ast}}  (y,z,q^{-1};\alpha_2)
$$
$$
 F_K (y,z,q;\alpha_2) \mapsto  F_{K^{\ast}}  (y^{-1} ,z^{-1} , q^{-1}; \alpha_2) =   F_{K^{\ast}}  (y,z,q^{-1};\alpha_1).
$$
Therefore, $F_{K^{\ast}}$ is defined by
$$
F_{K^{\ast}}  (y,z,q)  : =  F_K (y,z,q^{-1}) \in \intg + q^{-\Delta}\intg[q^{-1},q][[y/z ,\lb y/z \rb^{-1}]].
$$

\section{Surgery}

\subsection{Gluing}

\begin{figure}[t]
\begin{center}
\begin{tikzpicture}[scale=1]
\centering
\tikzstyle{every node}=[draw,shape=circle]

\draw (-1,0) node[circle,fill=white,inner sep=1pt,label=above:$w_1$](){}; 
\draw (-1,0)  -- (-2,0) node[circle,fill,inner sep=1pt,label=right:](){};
\draw (-2,0)  -- (-4,1);
\draw (-2,0) --  (-4,-1); 

\draw (-3,0.3) node[circle,fill,inner sep=1pt,label=right:](){};
\draw (-3,0) node[circle,fill,inner sep=1pt,label=right:](){};
\draw (-3,-0.3) node[circle,fill,inner sep=1pt,label=right:](){};

{$\bigcup$}

\draw (1,0) node[circle,fill=white,inner sep=1pt,label=above:$w_2$](){}; 
\draw (1,0)  -- (2,0) node[circle,fill,inner sep=1pt,label=right:](){};
\draw (2,0)  -- (4,1); 
\draw (2,0) --  (4,-1); 

\draw (3,0.3) node[circle,fill,inner sep=1pt,label=right:](){};
\draw (3,0) node[circle,fill,inner sep=1pt,label=right:](){};
\draw (3,-0.3) node[circle,fill,inner sep=1pt,label=right:](){};

\draw (4.7,0.1)  -- (5.2,0.1); 
\draw (4.7,-0.1)  -- (5.2,-0.1); 

\draw (9,0) node[circle,fill,inner sep=1pt,label=above:$w_1 + w_2 $](){}; 
\draw (9,0)  -- (10,0) node[circle,fill,inner sep=1pt,label=right:](){};
\draw (9,0)  -- (8,0) node[circle,fill,inner sep=1pt,label=left:](){};

\draw (10,0)  -- (12,1);
\draw (10,0) --  (12,-1); 

\draw (8,0)  -- (6,1); 
\draw (8,0) --  (6,-1); 

\draw (11,0.3) node[circle,fill,inner sep=1pt,label=right:](){};
\draw (11,0) node[circle,fill,inner sep=1pt,label=right:](){};
\draw (11,-0.3) node[circle,fill,inner sep=1pt,label=right:](){};

\draw (7,0.3) node[circle,fill,inner sep=1pt,label=right:](){};
\draw (7,0) node[circle,fill,inner sep=1pt,label=right:](){};
\draw (7,-0.3) node[circle,fill,inner sep=1pt,label=right:](){};

\end{tikzpicture}
\end{center}
\caption{Gluing of two plumbed knot complements results in a closed oriented plumbed 3-manifold. The union is along their boundary torus represented by the open vertices.}
\end{figure}
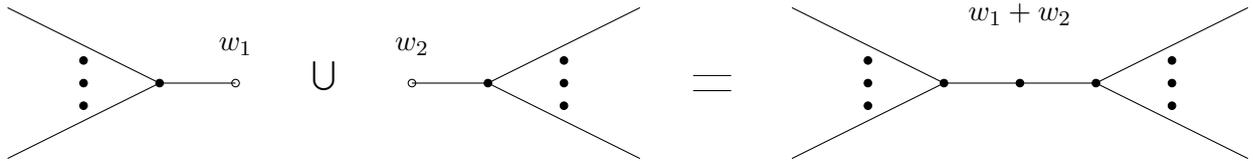

An important part of surgery is gluing. This procedure can produce a closed manifold when two manifolds with homeomorphic boundaries are attached. The resulting manifold depends on the details of the gluing.  In our setting, we have two knot complements. They can be glued along their common torus boundaries to obtain a closed oriented manifold. In the case of plumbed knot complements, the gluing of two distinguished vertices results in a closed oriented manifold as shown in Figure 6. We denote the two plumbed knot complements by $Y_1 (\Gamma_1 , v_{1, s}) $ and $Y_2 (\Gamma_2 , v_{2, s})$\footnote{We assume that the  knot complements are weakly negative definite.}. 
$$
Y_1 \cup_{T^2} Y_2 = Y = Y(\Gamma)
$$
where $\Gamma$ is obtained from $\Gamma_1$ and $\Gamma_2$ by attaching them as shown in Figure 6.\\
We further denote the adjacency matrices of $Y_1, Y_2$ and $Y$ by $B_1, B_2$ and $B$, respectively. The resulting $B$ is obtained by~\cite{GM}
$$
B = 
\begin{pmatrix}
\hat{B}_1 & | & \vdots & | & 0\\[10pt]
\hline
\ast \dots \ast & | & m_{1,s}+m_{2,s} & | & \ast \dots \ast\\[10pt]
\hline
0 & | & \vdots & | & \hat{B}_2\\[10pt]
\end{pmatrix}
$$
where $\hat{B}_1$ and $\hat{B}_2$ are adjacency matrices of $\Gamma_1 \backslash \lac v_{1,s}\rac$ and $\Gamma_2 \backslash \lac v_{2,s}\rac$, respectively.
\newline

To analyze how the relative $Spin^c$ structures behave under gluing, it is more useful to analyze from the viewpoint of cycles in $H_1 (Y_i) \simeq H^2 (Y_i, \ptl Y_i), i=1,2$. We have a surjective map from the Mayer-Vietoris sequence of $(Y_1, Y_2, Y)$:
$$
\lb H_1 (Y_1) \times H_1 (Y_1) \rb \oplus \lb H_1 (Y_2) \times H_1 (Y_2) \rb \rightarrow H_1 (Y) \times H_1 (Y)
$$
that is given by
$$
\lb \lsb (b^{(1)}_1,\cdots,b^{(1)}_s , c^{(1)}_1,\cdots,c^{(1)}_s ) \rsb , \lsb (b^{(2)}_1,\cdots,b^{(2)}_t , c^{+}_1,\cdots,c^{(2)}_t) \rsb \rb \mapsto
$$
$$
\lsb (b^{(1)}_1,\cdots,b^{(1)}_s + b^{(2)}_1 ,  b^{(2)}_2,\cdots b^{(2)}_t, c^{(1)}_1,\cdots,c^{(1)}_s +  c^{(2)}_1   , c^{(2)}_2 \cdots, c^{(2)}_t) \rsb
$$
We check the well definedness of the map (i.e. that it is independent of the choice of representatives of the relative $Spin^c$ structures). We pick 1-cycles from $H_1 (T^2) \times  H_1 (T^2)$. When gluing two boundaries, orientation of one of them is reversed, which results in the orientation reversal of a meridian $\mu_1 = - \mu_2$. Consequently, the actions by the meridians from the pair of 1-cycles
$$
b^{(1)}_s \mapsto b^{(1)}_s +1 ,\quad c^{(1)}_s \mapsto c^{(1)}_s +1 , \qquad b^{(2)}_1 \mapsto b^{(2)}_1 -1 ,\quad c^{(2)}_1 \mapsto c^{(2)}_1 - 1
$$
thus, the image does not change. In the case of longitudes, their actions are
$$
\vec{b}^{(1)}  \mapsto \vec{b}^{(1)} + B\vec{e}^{(1)}_s ,\quad \vec{c}^{(1)} \mapsto  \vec{c}^{(1)} + B\vec{e}^{(1)}_s , \qquad \vec{b}^{(2)} \mapsto \vec{b}^{(2)}  + B\vec{e}^{(2)}_1, \quad \vec{c}^{(2)} \mapsto  \vec{c}^{(2)} + B\vec{e}^{(2)}_1 
$$
Adding an element in the image of $B$ does change the resulting $Spin^c$ structures of $Y$. 
\newline

This allows us to write down the gluing formula for the super $\hat{Z}$.
\begin{equation}
\hat{Z}_{b,c}[Y;q] = (-1)^{\tau} q^{\chi} \sum_{n,m} \int \frac{dy}{i2\pi y}\frac{dz}{i2\pi z} \hat{Z}^{(\alpha_i)}_{b_1 ,c_1} (Y_1 ; y,z,n,m,q) \hat{Z}^{(\alpha_i)}_{b_2 ,c_2} (Y_2 ; y,z,n,m,q)
\end{equation}
where 
$$
\tau = \Pi(Y) - \Pi(Y_1) -\Pi(Y_2), \quad \chi = -( \vec{b}, B^{-1} \vec{c}) + ( \vec{b}_1 , B^{-1} \vec{c}_1) +  ( \vec{b}_2 , B^{-1} \vec{c}_2) \in \rational
$$
for any choice of chamber $\alpha_i , i=\pm$.

\subsection{TQFT properties}

We package the ingredients in the previous sections into the framework of topological quantum field theories (TQFT). This provides evidence for the existence of a 3-dimensional non-semisimple TQFT. By the axioms of $n$-dimensional TQFTs~\cite{At,Se}, to an $(n-1)$-dimensional manifold, a vector space~\footnote{The definition does not assume finite dimensionality. It is a consequence of the coevaluation map (finiteness principle).} over a field $\mathbb{F}$ is assigned,
$$
M^{n-1} \mapsto V_{\mathbb{F}}.
$$
To an $n$-dimensional manifold (bordism), a linear map between tensor products of vector spaces is assigned,
$$
M^n \mapsto f: \bigotimes\limits_{i}   V_i  \rightarrow \bigotimes\limits_{r}   V_r ,
$$
where $i$ and $r$ run over the incoming and outgoing boundaries of $M^n$, respectively~\footnote{In case of $i=r=0$, $M^n$ is a closed $n$-manifold, which is a bordism from an empty $(n-1)$-manifold $\phi^{n-1}$ to $\phi^{n-1}$; an element of $\mathbb{F}$ is assigned.}.\\
\indent In our 3-dimensional setting, a vector space $H_{T^2}$ is attached to the torus boundary $T^2$ of the knot complement $Y_K$ equipped with relative $Spin^c$ structures and $\mathbb{F}$ is the Novikov field. Specifically, we have
$$
\hat{Z}_{b,c} (Y_K ; y,z,n,m,q) = \sum_{\substack{ (i,j)\in \intg_{\geq 0}^2 \\ (i,j) \neq (0,0) }} b(n,m;i,j;q) \lb \frac{y^i}{z^j} - \frac{z^j}{y^i} \rb ,
$$
where the variables $n$ and $m$ are associated with the longitude of $T^2$ whereas $i$ and $j$ correspond to the meridian of $T^2$. 
$$
b(n,m;i,j;q)  = \sum_{w\in \rational }  b(n,m;i,j,w) q^w  \in \mathbb{F},\qquad  b(n,m;i,j,w) \in \rational ,
$$
where $\mathbb{F}$ consists of $q$-series such that the set $\Omega = \lac w | b(n,m;i,j,w)  \neq 0 \rac \subset \rational $ is bounded below. The super $\hat{Z}$ is a vector in $H_{T^2}$,
\begin{equation*}
\hat{Z}_{b,c} (Y_K ; y,z,n,m,q) \in H_{T^2}.
\end{equation*}
For a closed (oriented) 3-manifold $Y$ equipped with $Spin^c$ structures, we have
$$
\hat{Z}_{b,c}(Y;q) \in \mathbb{F}.
$$
Furthermore, there is an inner product (bilinear pairing) on $H_{T^2}$
$$
\left\langle b_1 | b_2 \right\rangle := \sum_{(n,m)} \sum_{(i,j)}  b_1 (n,m,i,j;q) b_2 (n,m,i,j;q) \in \mathbb{F}.
$$
The gluing formula (35) can be expressed via the inner product
$$
\hat{Z}_{b,c}(Y;q) = (-1)^{\tau} q^{\chi} \left\langle \hat{Z}_{b_1 , c_1}(Y_1) | R \hat{Z}_{b_2 , c_2}(Y_2) \right\rangle .
$$
where $R$ is the orientation reversal map for the meridian
$$
R : H_{T^2} \rightarrow H_{T^2},\quad (Rb)(m,n;i,j;q) = b(m,n;-i,-j;q) .
$$ 
So far, we have considered bordisms with one boundary component of genus 1. In order to arrive at the complete structure of the TQFT, for example, we have to consider bordisms with multiple number of boundary components of genus 1, as well as those of higher genus. We hope to investigate them in the future.

\subsection{The Dehn surgery formula}

We apply the results of the previous sections to derive the Dehn surgery formula for the super $F_K$. We first review Dehn surgery briefly.
\newline

\indent Let $Y$ be a closed oriented manifold and $K$ be a knot in $Y$. We carve out a tubular neighborhood of $K$, which is diffeomorphic to $S^1 \times D^2$. This yields a compact oriented manifold $Y_K$ with a torus boundary. Then we glue a solid torus $S^1 \times D^2$ into $Y_K$ along a slope $p/r \in \rational \cup \lac \infty \rac$ via a diffeomorphism. When gluing, a meridian of the solid torus is mapped to $p \mu + r \lambda$ on $\ptl Y_K =T^2$, where $\mu$ is a meridian and $\lambda$ is a longitude of $T^2$. This results in a closed oriented manifold $Y_{p/r}$.
$$
Y_{p/r} = Y_K \cup_{T^2} S^1 \times D^2 .
$$
In our setting, we have a plumbed knot in $Y=S^3$ and the surgery slope $p/r$ is specified by the solid torus in Section 4.2. Under Dehn surgery, we have the following relation between the super $F_K$ and $\hat{Z}$. 
\newline
\begin{Thm}
Let $Y_K$ be the complement of a knot $K$ in the 3-sphere $S^3$ and let $Y_{p/r}$ be the result of Dehn surgery along $K$ with slope $p/r \in \rational^{\ast}$. Assume that $Y_K$ and $Y_{p/r}$ are represented by negative definite plumbings. Then the invariants of $Y_{p/r}$ are given by
$$
\hat{Z}_{b,c}[Y_{p/r};q] = (-1)^{\tau} \mathcal{L}^{(\alpha_i ;\, p/r )}_{b,c} \lsb F^{(\alpha_i)}_K (y,z,q) \rsb,
$$
where the Laplace transform for the $\alpha_+ $ chamber is
$$
\mathcal{L}^{(\alpha_+ ;\, p/r )}_{b,c} : y^{\alpha}z^{\beta}q^{\gamma} \mapsto q^{\gamma}
\begin{cases}
\sum\limits_{r_s = r_{s,min}}^{\infty}  q^{\frac{\beta(r \alpha + \epsilon r_s )}{p}}, & \text{if}\quad r \alpha + \epsilon r_s + b \in p\intg,\, r\beta + c \in p\intg\\
\sum\limits_{w_s = w_{s,min} }^{\infty} q^{\frac{\alpha(r \beta - \epsilon w_s )}{p}}, & \text{if}\quad   r \beta - \epsilon w_s + c \in p\intg ,\, r\alpha + b \in p\intg\\
0, & \text{otherwise}
\end{cases}
$$
and the Laplace transform for the $\alpha_- $ chamber is 
$$
\mathcal{L}^{(\alpha_- ;\, p/r )}_{b,c} : y^{\alpha}z^{\beta}q^{\gamma} \mapsto -q^{\gamma}
\begin{cases}
\sum\limits_{w^{\prime}_s =  w^{\prime}_{s,min}  }^{\infty} q^{\frac{\beta(r \alpha - \epsilon w^{\prime}_s )}{p}}, & \text{if}\quad   r \alpha - \epsilon w^{\prime}_s + b \in p\intg ,\, r\beta + c \in p\intg\\
\sum\limits_{r^{\prime}_s =  r^{\prime}_{s,min} }^{\infty}  q^{\frac{\alpha(r \beta + \epsilon r^{\prime}_s )}{p}}, & \text{if}\quad r \beta + \epsilon r^{\prime}_s + c \in p\intg,\, r\alpha + b \in p\intg\\
0, & \text{otherwise}
\end{cases}
$$
where $r_{s,min},\, r^{\prime}_{s,min} \geq 1, \, w_{s,min},\, w^{\prime}_{s,min} \geq 0$ and $\epsilon = sign(p) (-1)^{\pi + 1}$.
\end{Thm}
We observe a qualitative difference between the above surgery formula and the $sl(2)$ surgery formula (39) in Appendix B. The latter transforms a term into a single term whereas the former converts a term into a series when contributing.
\newline

\noindent \textit{Proof}. We pick the $\alpha_+$ chamber and set the first $\hat{Z}$ to be $F_K$ and the second $\hat{Z}$ to be the solid torus (22) in (35). Expressing $F_K$ as
$$
F^{(\alpha_+)}_K (y,z,q) = \sum\limits_{\alpha , \beta , \gamma} C_{\alpha\beta\gamma} y^{\alpha} z^{\beta} q ^{\gamma}
$$
where $\alpha \in \intg_{\geq 0}, \beta \in \intg_{\leq 0}, \gamma \in \intg$ and $(\alpha , \beta) \neq (0,0)$. After substitution, the integrand of (35) becomes
\begin{equation}
F^{(\alpha_+)}_K \hat{Z}^{(\alpha_+)}_{b,c} = (-1)^{\pi} \sum\limits_{\alpha , \beta , \gamma} C_{\alpha\beta\gamma} y^{\alpha} z^{\beta} q ^{\gamma} 	 \lsb \sum_{ \Lambda_{b,c}^{-,0}} y_1^{r_1} z_1^{g_1} q^{\frac{g_1}{p} \lb r r_1  - \epsilon r_s \rb} +  \sum_{\Lambda_{b,c}^{0,+}} y_1^{d_1}z_1^{w_1} q^{ \frac{d_1}{p} \lb r w_1  + \epsilon w_s \rb} \rsb
\end{equation}
The integrations in (35) fix some of the summation indices to be
\begin{align}
r_1 & = -\alpha, & g_1 & = -\beta \nonumber\\
d_1 & = -\alpha, & w_1 & = -\beta
\end{align}
Substituting (37) into (36) and (25), we arrive at
\begin{align*}
r\alpha + \epsilon r_s +b & \in p\intg , & r\beta + c & \in p\intg \\
r\beta - \epsilon w_s + c & \in p\intg , & r\alpha + b & \in p\intg
\end{align*}
and the Laplace transform for the $\alpha_{+}$ chamber. Furthermore, we deduce from the above that 
$$
b,c \in \intg\, \text{mod}\, p.
$$
The Laplace transform for the $\alpha_{-}$ chamber can be derived in the same way.
\newline

\begin{Conj}
Let $K \subset S^3$ be a knot and $S^{3}_{p/r}(K)$ be the result of Dehn surgery on $K$. For any choice of a good chamber $\alpha_i$,
$$
\hat{Z}_{b,c}[S^{3}_{p/r}(K) ; q] = (-1)^{\tau} \mathcal{L}^{(\alpha_i ;\, p/r )}_{b,c} \lsb F^{(\alpha_i)}_K (y,z,q) \rsb,
$$
provided that the right hand side is well-defined.
\end{Conj}

\begin{Rmk} The above well-definedness condition can restrict the range of surgery slopes $p/r$; the range depends on the specific behaviors of $f_{m,n}(q)$ for $K$.
\end{Rmk}

\subsection{Examples}

\begin{figure}[t]
\begin{center}
\begin{tikzpicture}[scale=1]
\centering
\tikzstyle{every node}=[draw,shape=circle]

\draw (0,0) node[circle,fill,inner sep=1pt,label=above:$-1$](){}; 
\draw (0,0)  -- (1,1) node[circle,fill,inner sep=1pt,label=right:$-3$](){};
\draw (0,0)  -- (-1,1) node[circle,fill,inner sep=1pt,label=left:$-2$](){};

\draw (0,0) -- (0,-1) node[circle,fill,inner sep=1pt,label=right:$-7$](){};

\end{tikzpicture}
\quad 
\begin{tikzpicture}[scale=1]
\centering
\tikzstyle{every node}=[draw,shape=circle]

\draw (0,0) node[circle,fill,inner sep=1pt,label=above:$-1$](){}; 
\draw (0,0)  -- (1,1) node[circle,fill,inner sep=1pt,label=right:$-3$](){};
\draw (0,0)  -- (-1,1) node[circle,fill,inner sep=1pt,label=left:$-2$](){};

\draw (0,0) -- (0,-1) node[circle,fill,inner sep=1pt,label=right:$-7$](){};

\draw (0,0) -- (0,-2) node[circle,fill,inner sep=1pt,label=right:$-2$](){};

\end{tikzpicture}
\quad 
\begin{tikzpicture}[scale=1]
\centering
\tikzstyle{every node}=[draw,shape=circle]

\draw (0,0) node[circle,fill,inner sep=1pt,label=above:$-1$](){}; 
\draw (0,0)  -- (1,1) node[circle,fill,inner sep=1pt,label=right:$-3$](){};
\draw (0,0)  -- (-1,1) node[circle,fill,inner sep=1pt,label=left:$-2$](){};

\draw (0,0) -- (0,-1) node[circle,fill,inner sep=1pt,label=right:$-7$](){};
\draw (0,0) -- (0,-2) node[circle,fill,inner sep=1pt,label=right:$-2$](){};
\draw (0,0) -- (0,-3) node[circle,fill,inner sep=1pt,label=right:$-2$](){};

\end{tikzpicture}
\end{center}
\caption{The plumbing graphs of $\Sigma(2,3,7)$ (left), $\Sigma(2,3,13)$ (middle) and $\Sigma(2,3,19)$ (right).}
\end{figure}
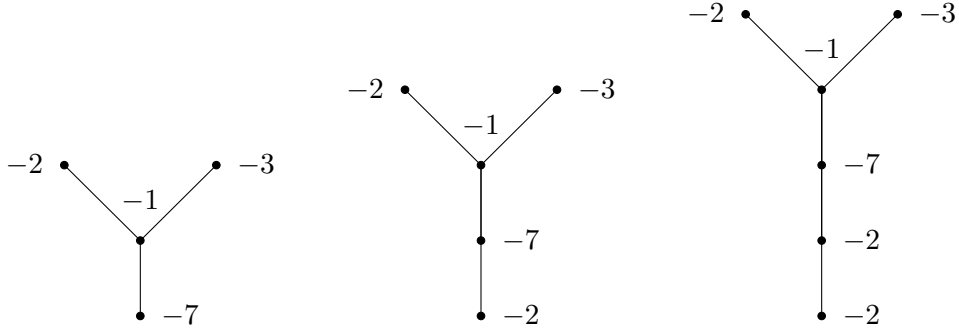

We apply Dehn surgery to torus knots. It is well known that this surgery produces a Seifert fibered manifold~\cite{M}. In particular, $-1/r$ surgery slopes yield the Brieskorn spheres, which are Seifert fibered integral homology spheres having three singular fibers.
$$
S^3_{-\frac{1}{r}}(T(s,t)) = \Sigma (s,t,rst+1),\qquad r \geq 1.
$$
Using Theorem 6.1 and Figure  4, we obtain the following $\hat{Z}$'s.
\newline

\noindent $S^3_{-1}(3_{1}^{r})=\Sigma(2,3,7)$
\begin{equation*}
\begin{split}
 \hat{Z} & \cong 2q^2+2 q^3+4 q^4+2 q^5+6 q^6+4 q^7+6 q^8+6 q^9+8 q^{10}+4 q^{11}+10 q^{12}+6 q^{13}+8 q^{14}+8 q^{15}\\
& +10 q^{16}+6q^{17}+12 q^{18}+ \cdots \\
\end{split}
\end{equation*}

\noindent $S^3_{-\frac{1}{2}}(3_{1}^{r})=\Sigma(2,3,13)$
\begin{equation*}
\begin{split}
 \hat{Z} & \cong 2q^2+2 q^3+4 q^4+2 q^5+6 q^6+2 q^7+6 q^8+4 q^9+6 q^{10}+2 q^{11}+10 q^{12}+4 q^{13}+6 q^{14}+8 q^{15}\\
& +10 q^{16}+4 q^{17}+10 q^{18} + \cdots  \\
\end{split}
\end{equation*}

\noindent $S^3_{-\frac{1}{3}}(3_{1}^{r})=\Sigma(2,3,19)$
\begin{equation*}
\begin{split}
 \hat{Z} & \cong 2q^2+2 q^3+4 q^4+2 q^5+6 q^6+2 q^7+6 q^8+4 q^9+6 q^{10}+2 q^{11}+10 q^{12}+4 q^{13}+6 q^{14}  +6 q^{15}\\
& +8 q^{16}+2 q^{17}+10 q^{18} + \cdots \\
\end{split}
\end{equation*}
The first example agrees with the result in \cite{FP}. The other examples coincide with the $\hat{Z}$ computed from Figure 7 using (14) and (17).
\begin{Rmk} The symbol $\cong$ denotes equality up to the additive constant ($\in \rational$) in (2).
\end{Rmk}

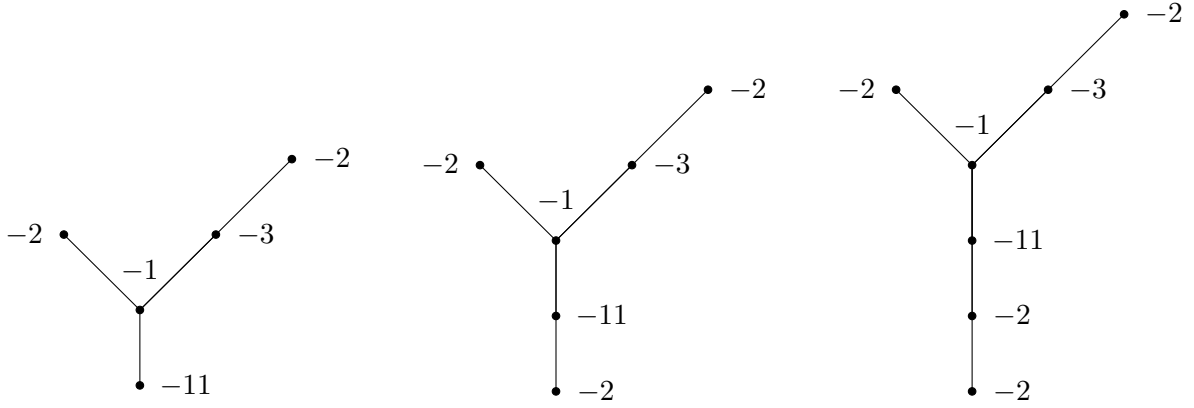
\begin{figure}[t]
\begin{center}
\begin{tikzpicture}[scale=1]
\centering
\tikzstyle{every node}=[draw,shape=circle]

\draw (0,0) node[circle,fill,inner sep=1pt,label=above:$-1$](-1){}; 
\draw (0,0)  -- (1,1) node[circle,fill,inner sep=1pt,label=right:$-3$](-3){};
\draw (0,0)  -- (2,2) node[circle,fill,inner sep=1pt,label=right:$-2$](){};
\draw (0,0)  -- (-1,1) node[circle,fill,inner sep=1pt,label=left:$-2$](){};

\draw (0,0) -- (0,-1) node[circle,fill,inner sep=1pt,label=right:$-11$](){};

\end{tikzpicture}
\quad 
\begin{tikzpicture}[scale=1]
\centering
\tikzstyle{every node}=[draw,shape=circle]

\draw (0,0) node[circle,fill,inner sep=1pt,label=above:$-1$](){};
\draw (0,0)  -- (1,1) node[circle,fill,inner sep=1pt,label=right:$-3$](){};
\draw (0,0)  -- (2,2) node[circle,fill,inner sep=1pt,label=right:$-2$](){};
\draw (0,0)  -- (-1,1) node[circle,fill,inner sep=1pt,label=left:$-2$](){};

\draw (0,0) -- (0,-1) node[circle,fill,inner sep=1pt,label=right:$-11$](){};
\draw (0,0) -- (0,-2) node[circle,fill,inner sep=1pt,label=right:$-2$](){};

\end{tikzpicture}
\quad 
\begin{tikzpicture}[scale=1]
\centering
\tikzstyle{every node}=[draw,shape=circle]

\draw (0,0) node[circle,fill,inner sep=1pt,label=above:$-1$](){};
\draw (0,0)  -- (1,1) node[circle,fill,inner sep=1pt,label=right:$-3$](){};
\draw (0,0)  -- (2,2) node[circle,fill,inner sep=1pt,label=right:$-2$](){};
\draw (0,0)  -- (-1,1) node[circle,fill,inner sep=1pt,label=left:$-2$](){};

\draw (0,0) -- (0,-1) node[circle,fill,inner sep=1pt,label=right:$-11$](){};
\draw (0,0) -- (0,-2) node[circle,fill,inner sep=1pt,label=right:$-2$](){};
\draw (0,0) -- (0,-3) node[circle,fill,inner sep=1pt,label=right:$-2$](){};

\end{tikzpicture}
\end{center}
\caption{The plumbing graphs of $\Sigma(2,5,11)$ (left), $\Sigma(2,5,21)$ (middle),  and $\Sigma(2,5,31)$ (right).}
\end{figure}

\noindent $S^3_{-1}(T(2,5))=\Sigma(2,5,11)$
\begin{equation*}
\begin{split}
 \hat{Z} & \cong 2 q^2+4 q^4+2 q^5+4 q^6+2 q^7+6 q^8+2 q^9+6 q^{10}+4 q^{11}+8 q^{12}+4 q^{13}+6 q^{14}+6 q^{15}\\
&+10 q^{16}+4 q^{17}+8q^{18}+ \cdots \\
\end{split}
\end{equation*}
\noindent $S^3_{-\frac{1}{2}}(T(2,5))=\Sigma(2,5,21)$
\begin{equation*}
\begin{split}
 \hat{Z} & \cong 2 q^2+4 q^4+2 q^5+4 q^6+2 q^7+6 q^8+2 q^9+6 q^{10}+2 q^{11}+8 q^{12}+2 q^{13}+6 q^{14}+4 q^{15}\\
& +8 q^{16}+2 q^{17}+8q^{18} +2 q^{19}+10 q^{20}+6 q^{21}+6 q^{22} +4 q^{23}+12 q^{24}+6 q^{25} +\cdots \\
\end{split}
\end{equation*}
\noindent $S^3_{-\frac{1}{3}}(T(2,5))=\Sigma(2,5,31)$
\begin{equation*}
\begin{split}
 \hat{Z} & \cong 2 q^2+4 q^4+2 q^5+4 q^6+2 q^7+6 q^8+2 q^9+6 q^{10}+2 q^{11}+8 q^{12} +2 q^{13}+6 q^{14}+4 q^{15} \\
& +8 q^{16}+2 q^{17}+8q^{18} +2 q^{19}+10 q^{20}+4 q^{21}+6 q^{22} +2 q^{23}+12 q^{24}+4 q^{25} +\cdots \\
\end{split}
\end{equation*}
\vspace{0.5cm}
The above results are in agreement with the $\hat{Z}$ computed from Figure 8 using (14) and (17) .  
\newline
\noindent We further verified the super $\hat{Z}$'s for $S^3_{-1/r}(K),\, K=T(3,4), T(3,5),T(3,7),\, r=1,2$ against the plumbing graph method (their super $F_K$'s are recorded in Appendix A). 
\newline

We next apply the surgery to the left handed trefoil.
\newline
\noindent $S^3_{-1}(3_{1}^{l})=\Sigma(2,3,5)$
\begin{equation*}
\begin{split}
 \hat{Z} &  \cong 2 q^2+2 q^3+4 q^4+4 q^5+6 q^6+4 q^7+8 q^8+6 q^9+8 q^{10}+6 q^{11}+10 q^{12}+6 q^{13}+10 q^{14}+8 q^{15}\\
& +10 q^{16}+6q^{17}+12 q^{18} + \cdots \\
\end{split}
\end{equation*}

\noindent $S^3_{-\frac{1}{2}}(3_{1}^{l})=\Sigma(2,3,11)$
\begin{equation*}
\begin{split}
 \hat{Z} & \cong 2 q^2+2 q^3+4 q^4+2 q^5+6 q^6+2 q^7+6 q^8+4 q^9+6 q^{10}+4 q^{11}+10 q^{12}+4 q^{13}+8 q^{14}+8 q^{15}\\
& +8 q^{16}+6q^{17}+10 q^{18}+ \cdots  \\
\end{split}
\end{equation*}

\noindent $S^3_{-\frac{1}{3}}(3_{1}^{l})=\Sigma(2,3,17)$
\begin{equation*}
\begin{split}
 \hat{Z} & \cong 2 q^2+2 q^3+4 q^4+2 q^5+6 q^6+2 q^7+6 q^8+4 q^9+6 q^{10}+2 q^{11}+10 q^{12}+2 q^{13}+6 q^{14}+6 q^{15}\\
& +8 q^{16}+4 q^{17}+10 q^{18} + \cdots \\
\end{split}
\end{equation*}
The first example agrees with the result in \cite{FP}. The other examples coincide with the results obtained from (14) and (17).
\newline

\noindent We next consider integer surgeries on the (0-framed) unknot.
$$
S^3_{-2}(\text{unknot})=L(2,1),
$$
using (19) and Theorem 6.1 yields
$$
\hat{Z} \cong 
\begin{cases}
2 q^2+4 q^4+4 q^6+6 q^8+4 q^{10}+8 q^{12}+4 q^{14}+8 q^{16}+6 q^{18}+8 q^{20}+4 q^{22}+12 q^{24} +\cdots,\\
q^{\frac{1}{2}}\lb 1+2 q+2 q^2+2 q^3+3 q^4+2 q^5+2 q^6+4 q^7+2 q^8+2 q^9+4 q^{10}+2 q^{11}+3 q^{12}+4 q^{13} +\cdots\rb\\
2q+2 q^2+4 q^3+2 q^4+4 q^5+4 q^6+4 q^7+2 q^8+6 q^9+4 q^{10}+4 q^{11}+4 q^{12}+4 q^{13}+4 q^{14}  +\cdots \\
2 q+2 q^2+4 q^3+2 q^4+4 q^5+4 q^6+4 q^7+2 q^8+6 q^9+4 q^{10}+4 q^{11}+4 q^{12}+4 q^{13}+4 q^{14}  +\cdots \\
\end{cases}
$$
This result agrees with that of \cite{FP}.
\vspace{0.5cm}
$$
S^3_{-3}(\text{unknot})=L(3,1),
$$
$$
\hat{Z} \cong 
\begin{cases}
2q^3+4 q^6+4 q^9+6 q^{12}+4 q^{15}+8 q^{18}+4 q^{21}+8 q^{24}+6 q^{27}+8 q^{30}+4 q^{33}+12 q^{36} +\cdots,\\
q^{\frac{1}{3}}\lb 2+4 q+4 q^2+4 q^3+4 q^4+6 q^5+4 q^6+4 q^7+4 q^8+8 q^9+4 q^{10}+4 q^{11}+4 q^{12}+8 q^{13} +\cdots\rb\\
q^{\frac{4}{3}}\lb 2+4 q^2+4 q^4+4 q^6+2 q^7+4 q^8+4 q^{10}+8 q^{12}+4 q^{14}+4 q^{16}+4 q^{17}+4 q^{18} +\cdots\rb\\
2 q+2 q^2+2 q^3+4 q^4+2 q^5+2 q^6+4 q^7+4 q^8+2 q^9+4 q^{10}+2 q^{11}+4 q^{12}+4 q^{13}+4 q^{14} +\cdots \\
2 q^2+2 q^4+2 q^5+2 q^6+4 q^8+4 q^{10}+2 q^{11}+2 q^{12}+4 q^{14}+2 q^{15}+4 q^{16}+2 q^{17} +2 q^{18}  +\cdots\\
q^{\frac{2}{3}}\lb 2+2 q+4 q^2+2 q^3+4 q^4+2 q^5+6 q^6+2 q^7+4 q^8+2 q^9+6 q^{10}+4 q^{11}+4 q^{12}+2 q^{13}  +\cdots\rb\\
\end{cases}
$$
This result coincides with that of \cite{FP}.
\newline

\noindent We consider integer surgeries on $T(2,3): S^3_{-p}(T(2,3))$.
\newline
$
p=2,\qquad  M \lb -1 \bigg| \frac{1}{2}, \frac{1}{3}, \frac{1}{8}\rb
$
$$
\hat{Z} \cong 
\begin{cases}
2 q^2+4 q^4+2 q^5+4 q^6+2 q^7+6 q^8+2 q^9+6 q^{10}+2 q^{11} +\cdots\\
2 q^{3/2} \lb 1+q+q^2+2 q^3+2 q^4+q^5+3 q^6+2 q^7+2 q^8+3 q^9+2 q^{10}+3 q^{11}+3 q^{12} +\cdots\rb\\
q+q^2+3 q^3+2 q^4+3 q^5+4 q^6+4 q^7+3 q^8+5 q^9+5 q^{10}+4 q^{11}  +\cdots \\
\end{cases}
$$
\newline
$
p=3,\qquad M \lb -1 \bigg| \frac{1}{2}, \frac{1}{3}, \frac{1}{9}\rb
$
$$
\hat{Z} \cong 
\begin{cases}
2 q^3+2 q^4+4 q^6+2 q^7+2 q^8+4 q^9+2 q^{10}+2 q^{11}+6 q^{12}+2 q^{13}+2 q^{14} +\cdots\\
q^{2/3} \left(1+q+3 q^2+2 q^3+3 q^4+2 q^5+5 q^6+2 q^7  +\cdots \right) \\
q+q^2+q^3+3 q^4+2 q^5+2 q^6+3 q^7+4 q^8+2 q^9+3 q^{10}+q^{11}+3 q^{12}  +\cdots\\
2 q^2+q^3+2 q^4+3 q^5+3 q^6+2 q^7+4 q^8+2 q^9+5 q^{10}+3 q^{11}   +\cdots\\
2 q^{4/3} \left(1+2 q^2+q^3+2 q^4+q^5+2 q^6+2 q^7+2 q^8+q^9  +\cdots\right)\\
2 q^{4/3} \left(1+q+q^2+q^3+2 q^4+q^5+2 q^6+q^7 +\cdots\right)\\
\end{cases}
$$
The results are in agreement with that of (14) and (17).
\begin{Rmk} In the above examples, the values of $(b,c)$ for $\hat{Z}$ in Theorem 6.1 are different from those of (14). An important point is that the application of Theorem 6.1 yields all the $q$-series in (14).
\end{Rmk}

\section{Open problems}

\begin{itemize}
	\item Finding a closed form formula for the super $F_K$ of all $T(s,t)$ that does not involve an algorithm would be valuable. Such a formula would provide an efficient way to obtain the super $F_K$ of these knots. Furthermore, torus knots are useful for a variety of purposes, as demonstrated by $F_K(x,q)$ associated with a Lie algebra. 
	
	\item more tractable problem than the one above is finding a general algorithm or a formula for $\epsilon_{m,n}$ functions for all torus knots. As we saw in Section 5.3, the $T(2,2l+1)$ family exhibits a pattern. Since the behavior of torus knots is uniform, we expect such an algorithm or a formula exists.
	
	\item Finding a super $\hat{Z}$ formula for positive definite plumbed manifolds is an open problem, and this approach appears to be challenging. An alternative route is via Dehn surgery. For this approach to be effective, a surgery formula that works for any positive surgery slope is essential. (cf. Remark 6.3). 
	
	\item To further develop TQFT properties, the super $F_K$ for higher genus surfaces ($g>1$) is necessary and valuable.
	
\end{itemize}

\appendix
\section*{Appendix}
\addcontentsline{toc}{section}{Appendix}


\section{Further examples}

We record information for $T(2,7)$ and other torus knots, applying the algorithm described in Section 5.3 to $T(2,7)$.

\begin{enumerate}
	\item  $g=\lac 1,3,5 \rac$ 
	 	\begin{align*}
			  r_{m}^{++} & = m - (14+2+7)\\
				r_{m}^{--} & = m - (14-2-7)\\
				r_{m}^{+-} & = m - (14+2-7)\\
				r_{m}^{-+} & = m - (14-2+7)\\			
			\end{align*}

		\item For $g=5$, $r_{m}^{++} \equiv 0, r_{m}^{--}\equiv 4$.
		
		\item $r_{n}^{++} \equiv 5,  r_{m}^{--}\equiv 4+5 =9$. Then $p=(5,9)$. These together with the above are conditions for $g=5$	
			
		\item For $g=3$, $r_{n}^{++} \equiv 5, r_{n}^{--}\equiv 9$ and $r_{m}^{++} \equiv 5-3=2, r_{m}^{--}\equiv 9-3=6$.
		
		\item For $g=1$, $r_{n}^{++} \equiv 5,  r_{n}^{--}\equiv 9$, $r_{m}^{++} \equiv 5-1=4, r_{m}^{--}\equiv 9-1=8$.
		
	  \item Another condition for $g=1$, $r_{m}^{++} \equiv 5, r_{m}^{--}\equiv 9$ and $r_{n}^{++} \equiv 6,  r_{n}^{--}\equiv 10$.
		
		\item Another condition for $g=3$, $r_{m}^{++} \equiv 5, r_{m}^{--}\equiv 9$ and $r_{n}^{++} \equiv 8,  r_{n}^{--}\equiv 12$.
		
		\item Another condition for $g=5$, $r_{m}^{++} \equiv 5, r_{m}^{--}\equiv 9$ and $r_{n}^{++} \equiv 10,  r_{n}^{--}\equiv 0$.
		
		\item For $g=5$, $r_{m}^{+-} \equiv 7, r_{m}^{-+}\equiv 11$.
		
		\item $r_{n}^{+-} \equiv 7+5 =12, r_{n}^{-+}\equiv 2$. Then $s=(12,2)$. These together with the above are conditions for $g=5$	
		
		\item For $g=3$, $r_{n}^{+-} \equiv 12, r_{n}^{-+}\equiv 2$ and $r_{m}^{+-} \equiv 12-3=9, r_{m}^{-+}\equiv 9+4=13$.
		
		\item For $g=1$, $r_{n}^{+-} \equiv 12, r_{n}^{-+}\equiv 2$ and $r_{m}^{+-} \equiv 12-1=11, r_{m}^{-+}\equiv 11+4=15\equiv 1$.
		
		\item Another condition for $g=1$, $r_{m}^{+-} \equiv 12, r_{m}^{-+}\equiv 2$ and $r_{n}^{+-} \equiv 13, r_{n}^{-+}\equiv 3$.
		
		\item Another condition for $g=3$, $r_{m}^{+-} \equiv 12, r_{m}^{-+}\equiv 2$ and $r_{n}^{+-} \equiv 15, r_{n}^{-+}\equiv 5$.
		
		\item Another condition for $g=5$, $r_{m}^{+-} \equiv 12, r_{m}^{-+}\equiv 2$ and $r_{n}^{+-} \equiv 3, r_{n}^{-+}\equiv 7$.
	
\end{enumerate}
We arrive at the following.
$$
\epsilon_{m,n}(T(2,7)) = \begin{cases}
+1, \, r^{+-}_m \equiv 7 \quad\& \quad  r^{+-}_n \equiv 12 \quad\& \quad   r^{-+}_m \equiv 11 \quad\& \quad r^{-+}_n \equiv 2 \quad \text{mod} \, 14\\
+1, \, r^{+-}_m \equiv 9 \quad\& \quad  r^{+-}_n \equiv 12 \quad\& \quad   r^{-+}_m \equiv 13 \quad\& \quad r^{-+}_n \equiv 2  \quad \text{mod} \, 14 \\
+1, \, r^{+-}_m \equiv 11 \quad\& \quad  r^{+-}_n \equiv 12 \quad\& \quad   r^{-+}_m \equiv 1 \quad\& \quad r^{-+}_n \equiv 2  \quad \text{mod} \, 14 \\
+1, \, r^{+-}_m \equiv 12 \quad\& \quad  r^{+-}_n \equiv 13 \quad\& \quad   r^{-+}_m \equiv 2 \quad\& \quad r^{-+}_n \equiv 3  \quad \text{mod} \, 14 \\
+1, \, r^{+-}_m \equiv 12 \quad\& \quad  r^{+-}_n \equiv 1 \quad\& \quad   r^{-+}_m \equiv 2 \quad\& \quad r^{-+}_n \equiv 5  \quad \text{mod} \, 14 \\
+1, \, r^{+-}_m \equiv 12 \quad\& \quad  r^{+-}_n \equiv 3 \quad\& \quad   r^{-+}_m \equiv 2 \quad\& \quad r^{-+}_n \equiv 7  \quad \text{mod} \, 14 \\
-1, \, r^{++}_m \equiv 0 \quad\& \quad  r^{++}_n \equiv 5 \quad\& \quad   r^{--}_m \equiv 4 \quad\& \quad  r^{--}_n \equiv 9  \quad \text{mod} \, 14 \\
-1, \, r^{++}_m \equiv 2 \quad\& \quad  r^{++}_n \equiv 5 \quad\& \quad   r^{--}_m \equiv 6 \quad\& \quad  r^{--}_n \equiv 9   \quad \text{mod} \,14 \\
-1, \, r^{++}_m \equiv 4 \quad\& \quad  r^{++}_n \equiv 5 \quad\& \quad   r^{--}_m \equiv 8 \quad\& \quad  r^{--}_n \equiv 9   \quad \text{mod} \,14 \\
-1, \, r^{++}_m \equiv 5 \quad\& \quad  r^{++}_n \equiv 6 \quad\& \quad   r^{--}_m \equiv 9 \quad\& \quad  r^{--}_n \equiv 10   \quad \text{mod} \,14 \\
-1, \, r^{++}_m \equiv 5 \quad\& \quad  r^{++}_n \equiv 8 \quad\& \quad   r^{--}_m \equiv 9 \quad\& \quad  r^{--}_n \equiv 12   \quad \text{mod} \,14 \\
-1, \, r^{++}_m \equiv 5 \quad\& \quad  r^{++}_n \equiv 10 \quad\& \quad   r^{--}_m \equiv 9 \quad\& \quad  r^{--}_n \equiv 0   \quad \text{mod} \,14 \\
0,  \, \text{otherwise}
\end{cases}
$$
where $r^{++}_m = m - (14+2+7), r^{--}_m = m - (14-2-7), r^{+-}_m = m - (14+2-7), r^{-+}_m = m - (14-2+7)$ and $r_n$'s can be obtained by replacing $m$ by $n$.
$$
g(m,n) = \begin{cases}
5, \, r^{+-}_m \equiv 7 \quad\& \quad  r^{+-}_n \equiv 12 \quad\& \quad   r^{-+}_m \equiv 11 \quad\& \quad r^{-+}_n \equiv 2 \quad \text{mod} \, 14\\
5, \, r^{+-}_m \equiv 12 \quad\& \quad  r^{+-}_n \equiv 3 \quad\& \quad   r^{-+}_m \equiv 2 \quad\& \quad r^{-+}_n \equiv 7 \quad \text{mod} \, 14\\
5, \, r^{++}_m \equiv 0 \quad\& \quad  r^{++}_n \equiv 5 \quad\& \quad   r^{--}_m \equiv 4 \quad\& \quad  r^{--}_n \equiv 9  \quad \text{mod} \, 14 \\
5, \, r^{++}_m \equiv 5 \quad\& \quad  r^{++}_n \equiv 10 \quad\& \quad   r^{--}_m \equiv 9 \quad\& \quad  r^{--}_n \equiv 0  \quad \text{mod} \, 14 \\
3, \, r^{+-}_m \equiv 9 \quad\& \quad  r^{+-}_n \equiv 12 \quad\& \quad   r^{-+}_m \equiv 13 \quad\& \quad r^{-+}_n \equiv 2 \quad \text{mod} \, 14\\
3, \, r^{+-}_m \equiv 12 \quad\& \quad  r^{+-}_n \equiv 1 \quad\& \quad   r^{-+}_m \equiv 2 \quad\& \quad r^{-+}_n \equiv 5  \quad \text{mod} \, 14 \\
3, \, r^{++}_m \equiv 2 \quad\& \quad  r^{++}_n \equiv 5 \quad\& \quad   r^{--}_m \equiv 6 \quad\& \quad  r^{--}_n \equiv 9  \quad \text{mod} \, 14 \\
3, \, r^{++}_m \equiv 5 \quad\& \quad  r^{++}_n \equiv 8 \quad\& \quad   r^{--}_m \equiv 9 \quad\& \quad  r^{--}_n \equiv 12   \quad \text{mod} \,14 \\
1, \, r^{+-}_m \equiv 11 \quad\& \quad  r^{+-}_n \equiv 12 \quad\& \quad   r^{-+}_m \equiv 1 \quad\& \quad r^{-+}_n \equiv 2  \quad \text{mod} \, 14 \\
1, \,  r^{+-}_m \equiv 12 \quad\& \quad  r^{+-}_n \equiv 13 \quad\& \quad   r^{-+}_m \equiv 2 \quad\& \quad r^{-+}_n \equiv 3  \quad \text{mod} \, 14 \\
1, \, r^{++}_m \equiv 5 \quad\& \quad  r^{++}_n \equiv 6 \quad\& \quad   r^{--}_m \equiv 9 \quad\& \quad  r^{--}_n \equiv 10   \quad \text{mod} \,14 \\
1, \, r^{++}_m \equiv 4 \quad\& \quad  r^{++}_n \equiv 5 \quad\& \quad   r^{--}_m \equiv 8 \quad\& \quad  r^{--}_n \equiv 9   \quad \text{mod} \,14 \\
0, \, \text{otherwise}
\end{cases}
$$

$$
F_{T(3,4)}(y,z,q) = 1 + \sum_{\substack{i=3 \\ i\neq 5}}^{\infty} \lb y^i + \frac{1}{z^i} \rb  - \sum_{\substack{i=3 \\ i\neq 5}}^{\infty} \lb \frac{1}{y^i} + z^i \rb 
+q \left(\frac{y^3}{z^4}+\frac{y^4}{z^3}-\frac{z^3}{y^4}-\frac{z^4}{y^3}\right)
$$
$$
+q^2 \left(\frac{y^3}{z^8}+\frac{y^4}{z^6}+\frac{y^6}{z^4}+\frac{y^8}{z^3}-\frac{z^3}{y^8}-\frac{z^4}{y^6}-\frac{z^6}{y^4}-\frac{z^8}{y^3}\right)+q^3 \left(\frac{y^4}{z^9}+\frac{y^9}{z^4}-\frac{z^4}{y^9}-\frac{z^9}{y^4}\right)+q^4\left(\frac{y^6}{z^8}+\frac{y^8}{z^6}-\frac{z^6}{y^8}-\frac{z^8}{y^6}\right)
$$
$$
+q^6 \left(\frac{y^8}{z^9}+\frac{y^9}{z^8}-\frac{z^8}{y^9}-\frac{z^9}{y^8}\right)+q^7 \left(-\frac{y^7}{z^{12}}-\frac{y^{12}}{z^7}+\frac{z^7}{y^{12}}+\frac{z^{12}}{y^7}\right)+q^{10}\left(-\frac{y^{10}}{z^{12}}-\frac{y^{12}}{z^{10}}+\frac{z^{10}}{y^{12}}+\frac{z^{12}}{y^{10}}\right)
$$
$$
+q^{11} \left(-\frac{y^{11}}{z^{12}}-\frac{y^{12}}{z^{11}}+\frac{z^{11}}{y^{12}}+\frac{z^{12}}{y^{11}}\right)+q^{13}\left(-\frac{y^{12}}{z^{13}}-\frac{y^{13}}{z^{12}}+\frac{z^{12}}{y^{13}}+\frac{z^{13}}{y^{12}}\right)+q^{14}\left(-\frac{y^{12}}{z^{14}}-\frac{y^{14}}{z^{12}}+\frac{z^{12}}{y^{14}}+\frac{z^{14}}{y^{12}}\right)
$$
$$
+q^{17}\left(-\frac{y^{12}}{z^{17}}-\frac{y^{17}}{z^{12}}+\frac{z^{12}}{y^{17}}+\frac{z^{17}}{y^{12}}\right) + \cdots
$$

$$
F_{T(3,5)}(y,z,q) = 1 + \sum_{\substack{i=3 \\ i\neq 4,7}}^{\infty} \lb y^i + \frac{1}{z^i} \rb  - \sum_{\substack{i=3 \\ i\neq 4,7}}^{\infty} \lb \frac{1}{y^i} + z^i \rb + q \left(\frac{y^3}{z^5}+\frac{y^5}{z^3}-\frac{z^3}{y^5}-\frac{z^5}{y^3}\right)
$$
$$
+q^2\left(\frac{y^3}{z^{10}}+\frac{y^5}{z^6}+\frac{y^6}{z^5}+\frac{y^{10}}{z^3}-\frac{z^3}{y^{10}}-\frac{z^5}{y^6}-\frac{z^6}{y^5}-\frac{z^{10}}{y^3}\right)+q^3 \left(\frac{y^5}{z^9}+\frac{y^9}{z^5}-\frac{z^5}{y^9}-\frac{z^9}{y^5}\right)
$$
$$
+q^4\left(\frac{y^5}{z^{12}}+\frac{y^6}{z^{10}}+\frac{y^{10}}{z^6}+\frac{y^{12}}{z^5}-\frac{z^5}{y^{12}}-\frac{z^6}{y^{10}}-\frac{z^{10}}{y^6}-\frac{z^{12}}{y^5}\right)+q^6\left(\frac{y^9}{z^{10}}+\frac{y^{10}}{z^9}-\frac{z^9}{y^{10}}-\frac{z^{10}}{y^9}\right)
$$
$$
+q^8 \left(-\frac{y^8}{z^{15}}+\frac{y^{10}}{z^{12}}+\frac{y^{12}}{z^{10}}-\frac{y^{15}}{z^8}+\frac{z^8}{y^{15}}-\frac{z^{10}}{y^{12}}-\frac{z^{12}}{y^{10}}+\frac{z^{15}}{y^8}\right)+q^{11} \left(-\frac{y^{11}}{z^{15}}-\frac{y^{15}}{z^{11}}+\frac{z^{11}}{y^{15}}+\frac{z^{15}}{y^{11}}\right)
$$
$$
+q^{13} \left(-\frac{y^{13}}{z^{15}}-\frac{y^{15}}{z^{13}}+\frac{z^{13}}{y^{15}}+\frac{z^{15}}{y^{13}}\right)+q^{14} \left(-\frac{y^{14}}{z^{15}}-\frac{y^{15}}{z^{14}}+\frac{z^{14}}{y^{15}}+\frac{z^{15}}{y^{14}}\right) +q^{16} \left(-\frac{y^{15}}{z^{16}}-\frac{y^{16}}{z^{15}}+\frac{z^{15}}{y^{16}}+\frac{z^{16}}{y^{15}}\right)
$$
$$
+q^{17}\left(-\frac{y^{15}}{z^{17}}-\frac{y^{17}}{z^{15}}+\frac{z^{15}}{y^{17}}+\frac{z^{17}}{y^{15}}\right) +q^{19}\left(-\frac{y^{15}}{z^{19}}-\frac{y^{19}}{z^{15}}+\frac{z^{15}}{y^{19}}+\frac{z^{19}}{y^{15}}\right)+\cdots
$$

$$
F_{T(3,7)}(y,z,q) = 1 + \sum_{\substack{i=3 \\ i\neq 4,5,8,11}}^{\infty} \lb y^i + \frac{1}{z^i} \rb  - \sum_{\substack{i=3 \\ i\neq 4,5,8,11}}^{\infty} \lb \frac{1}{y^i} + z^i \rb 
+ q \left(\frac{y^3}{z^7}+\frac{y^7}{z^3}-\frac{z^3}{y^7}-\frac{z^7}{y^3}\right)
$$
$$
+q^2\left(\frac{y^3}{z^{14}}+\frac{y^6}{z^7}+\frac{y^7}{z^6}+\frac{y^{14}}{z^3}-\frac{z^3}{y^{14}}-\frac{z^6}{y^7}-\frac{z^7}{y^6}-\frac{z^{14}}{y^3}\right)+q^3
   \left(\frac{y^7}{z^9}+\frac{y^9}{z^7}-\frac{z^7}{y^9}-\frac{z^9}{y^7}\right)
$$
$$
+q^4 \left(\frac{y^6}{z^{14}}+\frac{y^7}{z^{12}}+\frac{y^{12}}{z^7}+\frac{y^{14}}{z^6}-\frac{z^6}{y^{14}}-\frac{z^7}{y^{12}}-\frac{z^{12}}{y^7}-\frac{z^{14}}{y^6}\right) +q^5 \left(\frac{y^7}{z^{15}}+\frac{y^{15}}{z^7}-\frac{z^7}{y^{15}}-\frac{z^{15}}{y^7}\right)
$$
$$
+q^6\left(\frac{y^7}{z^{18}}+\frac{y^9}{z^{14}}+\frac{y^{14}}{z^9}+\frac{y^{18}}{z^7}-\frac{z^7}{y^{18}}-\frac{z^9}{y^{14}}-\frac{z^{14}}{y^9}-\frac{z^{18}}{y^7}\right)+q^8 \left(\frac{y^{12}}{z^{14}}+\frac{y^{14}}{z^{12}}-\frac{z^{12}}{y^{14}}-\frac{z^{14}}{y^{12}}\right)
$$
$$
+q^{10} \left(-\frac{y^{10}}{z^{21}}+\frac{y^{14}}{z^{15}}+\frac{y^{15}}{z^{14}}-\frac{y^{21}}{z^{10}}+\frac{z^{10}}{y^{21}}-\frac{z^{14}}{y^{15}}-\frac{z^{15}}{y^{14}}+\frac{z^{21}}{y^{10}}\right)+q^{12} \left(\frac{y^{14}}{z^{18}}+\frac{y^{18}}{z^{14}}-\frac{z^{14}}{y^{18}}-\frac{z^{18}}{y^{14}}\right)
$$
$$
+q^{13}\left(-\frac{y^{13}}{z^{21}}-\frac{y^{21}}{z^{13}}+\frac{z^{13}}{y^{21}}+\frac{z^{21}}{y^{13}}\right) +q^{16} \left(-\frac{y^{16}}{z^{21}}-\frac{y^{21}}{z^{16}}+\frac{z^{16}}{y^{21}}+\frac{z^{21}}{y^{16}}\right) + \cdots
$$

\section{GM series}

We summarize the series invariant $F_K$ associated with a Lie algebra $sl(2)$ in \cite{GM} (overall $q$ factors and constants are suppressed).
\begin{equation}
F_K (x,q) = \sum_{\substack {m =1 \\ \text{odd} } }^{\infty} f_m (q) \lb x^{m/2} - x^{-m/2} \rb \quad \in q^{\Delta} \intg [q^{-1}, q]][[x^{1/2}, x^{-1/2}]]
\end{equation}
The Dehn surgery formula is given by
\begin{equation}
\hat{Z}_b (Y;q) \cong \mathcal{L}^{(b)}_{p/r} \lsb \lb x^{\frac{1}{2r}} - x^{-\frac{1}{2r}} \rb F_K (x,q) \rsb
\end{equation}
where
$$
\mathcal{L}^{(b)}_{p/r} :  x^{u} q^v  \mapsto
\begin{cases}
q^{-u^2 r/p} q^v , & \text{if}\quad ru-b \in p\intg \\
0,                & \text{otherwise}.
\end{cases}
$$
In case of the torus knots $K=T(s,t) \subset S^3$:
\begin{equation}
F_K (x,q) \cong  \sum_{\substack {m =1 \\ \text{odd} } }^{\infty} \epsilon_m  q^{\frac{m^2}{4st}}  \lb x^{m/2} - x^{-m/2} \rb
$$
where
$$
\epsilon_m =
\begin{cases}
+1 , & \text{if}\quad m \equiv st + s+t \quad\text{or}\quad st -s-t\quad  mod\, 2st\\
-1,  &  \text{if}\quad m \equiv st + s-t \quad\text{or}\quad st -s+t\quad  mod\, 2st\\
0,                 & \text{otherwise}.
\end{cases}
\end{equation}

\section{Supergroup Chern-Simons theory}

\begin{figure}
\centering

  \includegraphics[scale=1.2]{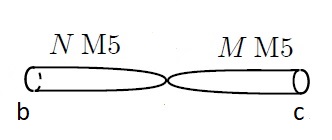}
	\qquad
%
  \includegraphics[scale=0.7]{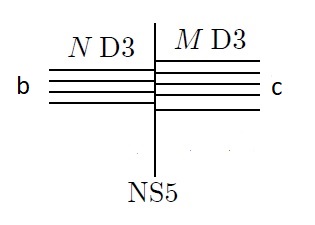}
\caption{The cigars of the Taub-NUT space of 11-dimensional spacetime that are wrapped by the branes (left). The brane system of the Type IIB theory (right). The labels $b$ and $c$ are the asymptotic boundary conditions taking values in  $H_1 (M^3 ;\intg)^N \times H_1 (M^3 ;\intg)^M$ for $U(N|M)$.}
\end{figure}

We review the physical aspects of $\hat{Z}_{b,c}$ including Chern-Simons theory on $Y$ associated with a Lie supergroup $U(N|M)$ in \cite{FP, MW} (see also \cite{V}).
\newline

We begin with a brane system in a 11-dimensional spacetime (ST) in M-theory. We take the 10d spatial geometry to be a cotangent bundle of a 3-manifold $M^3=Y$ and the 4-dimensional Taub-NUT (TN) space. The former is assumed to be a rational homology sphere. The latter looks like two cigars whose tips are joined at an origin. Away from the tip, the geometry looks like $S^{1}_M \times \real^3$, where the circle is taken to be the M-theory circle. Near tip geometry looks like $\complex^2 \cong \real^4$.
\begin{center}
\begin{tabular}{c c c c c c c c}
11D ST      &   $S^{1}_{\text{t}}$ &  $\times$ &  $T^{\ast}M^3$ & $\times$ & Taub           & $-$       & NUT\\ 
M M5 branes & $S^{1}_{\text{t}}$   & $\times$  & $M^3$          & $\times$ & $\complex$     & $\times $ & $\lac 0 \rac$\\  
N M5 branes & $S^{1}_{\text{t}}$   & $\times$  & $M^3$          & $\times$ & $\lac 0 \rac$  & $\times$  &  $\complex$
\end{tabular}
\end{center}
where $S^{1}_{\text{t}}$ is a time circle. The two stacks of M5 branes wrap the indicated parts of the spacetime as shown in Figure 9. The copies of $\complex $ are part of the TN space. This spacetime geometry has symmetries from the TN space, $U(1)_q \times U(1)_R$~\footnote{If $M^3$ has a circle fiber, for example, a Seifert fibered manifold, then an extra symmetry group $U(1)$ exists.}. We next shrink $S^{1}_M$ to reduce to 10 dimensional spacetime. This process lands us in type IIA string theory and the brane system becomes
\begin{center}
\begin{tabular}{c c c c c c}
Type IIA 10D ST  & $S^{1}_{\text{t}}$ &  $\times $ & $T^{\ast}M^3 $ & $\times $ & $\real^3$\\
1 D6 brane       & $S^{1}_{\text{t}}$ &  $\times $ & $T^{\ast}M^3 $ & $\times $ & $\lac 0 \rac$\\
M D4 branes      & $S^{1}_{\text{t}}$ &  $\times $ & $M^3$          & $\times $ & $\real_{+}$\\  
N D4 branes      & $S^{1}_{\text{t}}$ &  $\times $ & $M^3$          & $\times $ & $\real_{-}$ 
\end{tabular}
\end{center}
The M5 branes are transformed into the D4 branes. The D6 brane appears as a consequence of the Taub-NUT space. We apply T-duality along $S^{1}_{\text{t}}$ to pass to type IIB. Then we apply S-duality. We arrive at the following final brane system shown in Figure 9.
\begin{center}
\begin{tabular}{c c c c c c}
Type IIB 10D ST  & $S^{1}$     & $\times$ & $T^{\ast}M^3 $ & $\times $  & $\real^3$\\
1 NS5 brane      & $\text{pt}$ & $\times$ & $T^{\ast}M^3 $ & $\times $  & $\lac 0 \rac$\\
M D3 branes      & $\text{pt}$ & $\times$ & $M^3 $         & $\times $  & $\real_{+}$\\  
N D3 branes      & $\text{pt}$ & $\times$ & $M^3$          & $\times $  & $\real_{-}$ 
\end{tabular}
\end{center}
The S-duality maps D5 brane to NS5 brane. The former was obtained from the above D6 brane. On the stack of the D3 branes, its worldvolume theory is $4d\quad \mathcal{N}=4$ super Yang-Mills with gauge groups $U(M)$ whereas the theory on the other brane stack has gauge group $U(N)$.
\newline

We next apply the (GL) topological twist along $M^3$ of $T^{\ast} M^3$ to the above super Yang-Mills theories~\cite{KW}. This results in a cohomological quantum field theory that is a coupled 4d-3d system across the NS5 brane. The cohomological sector of the theory is the Chern-Simons theory based on $U(N|M)$ (up to $Q$-exact terms). Its action functional is the supergroup Chern-Simons theory (up to certain exact terms). Furthermore, analogous to the Chern-Simons level parameter in case of a Lie group $SU(N)$, $U(N|M)$ Chern-Simons theory carries a parameter $K$, which comes from the complexified gauge coupling constant $\tau$ of the super Yang-Mills theory, which in turn comes from the complexified string coupling constant.
$$
\tau = K cos(\theta) e^{i \theta} \in H^{+},
$$
where $\theta$ is the vaccum angle and $H^{+}$ the upper half complex plane ($Im\, \tau > 0$). The action functional of $U(M|N)$ Chern-Simons theory on $M^3$ at level $K$ is
$$
CS(\mathcal{A}) = \frac{iK}{4\pi} \int_{M^3} Str \lb \mathcal{A}  d \mathcal{A} + \frac{2}{3} \mathcal{A}^3 \rb + \lac Q, \cdots \rac, 
$$
where $\mathcal{A} = \mathcal{A}_b + \mathcal{A}_f, \mathcal{A}_b$ is the complexified gauge connection of $A$ and $\mathcal{A}_f$ is a fermion field. Str denotes the supertrace.
\newline

\noindent The existence of the super $\hat{Z}_{b,c}$ can be predicted from 11 dimensions. Specifically, the presence of the cigars in Figure 9, in particular their geometry away from the tips, requires imposing (asymptotic) boundary conditions $(b,c) \in H_1 (M^3 ;\intg)^N \times H_1 (M^3 ;\intg)^M$. The partition function over the BPS sector of the Hilbert space of the brane system is
$$
\hat{Z}^{gl(N|M)}_{b,c} [M^3 ; q] : = Tr_{H_{b,c}} (-1)^F q^{L_0}.
$$
where $F$ is fermion number operator and $L_0 $ is the generator of $U(1)_q$.

URL- https://sites.google.com/view/john-chae/home

\end{document}